\documentclass{amsart}
\usepackage{amssymb, amscd}
\usepackage{enumerate}
\usepackage{pstricks}

\newcommand\datver[1]{\def\datverp
 {\par\boxed{\boxed{\text{Version: #1; Run: \today}}}}}
\datver{3.1; Revised: March 22, 2006 by B.}
\newcommand{\ie}{{i.\thinspace e., }}

\newcommand{\CI}{\mathcal{C}^\infty}
\newcommand{\CIc}{{\mathcal C}^{\infty}_{\text{c}}}
\newcommand{\pa}{\partial}
\newcommand\<{\langle}
\renewcommand\>{\rangle}

\newcommand{\maD}{\mathcal D}

\newcommand{\maP}{\mathcal P}

\newcommand{\maO}{\mathcal O}

\newcommand{\maX}{\mathcal X}

\newcommand{\CC}{\mathbb C}
\newcommand{\NN}{\mathbb N}
\newcommand{\PP}{\mathbb P}
\newcommand{\RR}{\mathbb R}

\newcommand{\ZZ}{\mathbb Z}

\let\ep\epsilon
\newcommand{\patop}{\pa_{\rm top}}

\newcommand\comment[1]{}
\newcommand\Wkp{W^{k,p}}
\newcommand\WkpM{W^{k,p}(M_0)}
\newcommand\WkpO{W^{k,p}(\Omega_0)}
\newcommand\Wkinf{W^{k, \infty}}

\newcommand\smb{\nu^A}
\newcommand{\VV}{\mathcal V}
\newcommand{\WW}{\mathcal W}
\newcommand{\maK}{\mathcal K}

\newcommand{\Diff}[1]{{\rm Diff}(#1)}
\newcommand{\DiffV}[1]{{\rm Diff}^{#1}_{\VV}}

\newcommand{\vol}{\operatorname{vol}}
\newcommand{\supp}{\operatorname{supp}}
\newcommand\Kond[2]{{\mathcal K}^{#1}_{#2}}

\def\II{{\rm II}}

\newtheorem{theorem}{Theorem}[section]
\newtheorem{proposition}[theorem]{Proposition}
\newtheorem{corollary}[theorem]{Corollary}

\newtheorem{lemma}[theorem]{Lemma}

\theoremstyle{definition}
\newtheorem{definition}[theorem]{Definition}
\theoremstyle{remark}
\newtheorem{remark}[theorem]{Remark}
\newtheorem{example}[theorem]{Example}
\newtheorem{examples}[theorem]{Examples}

\author[B. Ammann]{Bernd Ammann} \address{Bernd Ammann,
Institut \'Elie Cartan, Universit\'e Henri Poincar\'e Nancy 1, B.P. 239,
 54506 Vandoeuvre-Les-Nancy, France,
{\rm http://www.differentialgeometrie.de/ammann}}
\email{bernd.ammann@gmx.net}

\author[A. Ionescu]{Alexandru D. Ionescu} \address{Alexandru D. Ionescu,
University of Wisconsin, Department of Mathematics, Madison, WI 53705, USA}
\email{ionescu@math.wisc.edu}

\author[V. Nistor]{Victor Nistor} \address{Victor Nistor, Pennsylvania State
University, Math. Dept., University Park, PA 16802, USA}
\email{nistor@math.psu.edu}

\thanks{Ionescu was supported in part by an NSF grant, by an Alfred
P. Sloan research fellowship, and by the Packard Foundation. 
Nistor was partially supported by the
NSF Grant DMS 0200808.}

\begin{document}

\date{June 21, 2006}

\title[Sobolev spaces]{Sobolev spaces on Lie manifolds and
regularity for polyhedral domains}

\begin{abstract}{We study some basic analytic questions
related to differential operators on Lie manifolds, which are
manifolds whose large scale geometry can be described by a a Lie
algebra of vector fields on a compactification. We extend to Lie
manifolds several classical results on Sobolev spaces, elliptic
regularity, and mapping properties of pseudodifferential
operators. A tubular neighborhood theorem for Lie submanifolds
allows us also to extend to regular open subsets of Lie manifolds
the classical results on traces of functions in suitable Sobolev
spaces. Our main application is a regularity result on polyhedral
domains $\PP \subset \RR^3$ using the weighted Sobolev spaces
$\Kond{m}{a}(\PP)$. In particular, we show that there is no loss
of $\Kond{m}{a}$--regularity for solutions of strongly elliptic
systems with smooth coefficients. For the proof, we identify
$\Kond{m}{a}(\PP)$ with the Sobolev spaces on $\PP$ associated to
the metric $r_{\PP}^{-2} g_E$, where $g_E$ is the Euclidean metric
and $r_{\PP}(x)$ is a smoothing of the Euclidean distance from $x$
to the set of singular points of $\PP$. A~suitable
compactification of the interior of $\PP$ then becomes a regular
open subset of a Lie manifold. We also obtain the well-posedness
of a non-standard boundary value problem on a smooth, bounded
domain with boundary $\maO \subset \RR^n$ using weighted Sobolev
spaces, where the weight is the distance to the boundary.}
\end{abstract}

\maketitle \tableofcontents


\section*{Introduction}

We study some basic analytic questions on non-compact
manifolds. In order to obtain stronger results, we restrict
ourselves to ``Lie manifolds,'' a class of manifolds whose large
scale geometry is determined by a compactification to a manifold
with corners and a Lie algebra of vector fields on this
compactification (Definition~\ref{def.unif.str}). One of the
motivations for studying Lie manifolds is the loss of (classical
Sobolev) regularity of solutions of elliptic equations on
non-smooth domains. To explain this loss of regularity, let us
recall first that the Poisson problem
\begin{equation}\label{eq.BVP}
    \Delta u = f \in H^{m-1}(\Omega), \quad m \in \NN\cup \{0\},\
        \Omega \subset
    \RR^n \text{ bounded,}
\end{equation}
has a unique solution $u \in H^{m+1}(\Omega)$, $u = 0$ on $\pa
\Omega$, provided that {\em $\pa \Omega$ is smooth}.  In particular,
$u$ will be smooth up to the boundary if $\pa \Omega$ and $f$ are
smooth (in the following, when dealing with functions defined on an
open set, by ``smooth,'' we shall mean ``smooth up to the
boundary''). See the books of Evans \cite{Evans}, or Taylor
\cite{Taylor1} for a proof of this basic well-posedness result.

This well-posedness result is especially useful in practice for the
numerical approximation of the solution $u$ of Equation \eqref{eq.BVP}
\cite{BabuAziz}. However, in practice, it is only rarely the case that
$\Omega$ is smooth. The lack of smoothness of the domains interesting
in applications has motivated important work on Lipschitz domains, see
for instance \cite{JerisonKenig, MitreaTaylor} or
\cite{Verchota}. These papers have extended to Lipschitz domains some
of the classical results on the Poisson problem on smooth, bounded
domains, using the classical Sobolev spaces
\begin{equation*}
    H^m(\Omega) := \{u,\ \pa^{\alpha} u \in L^2(\Omega),
    \ |\alpha| \le m \}.
\end{equation*}

%
It turns out that, if $\pa \Omega$ is {\em not} smooth, then the
smoothness of~$f$ on $\overline{\Omega}$ (\ie up to the boundary) does
not imply that the solution~$u$ of Equation~\eqref{eq.BVP} is smooth
as well on $\overline{\Omega}$. This is the {\em loss of regularity}
for elliptic problems on non-smooth domains mentioned above.

The loss of regularity can be avoided, however, by a conformal
blowup of the singular points. This conformal blowup replaces a
neighborhood of each connected component of the set of singular
boundary points by a complete, but non-compact end. (Here
``complete'' means complete as a metric space, not geodesically
complete.) It can be proved then that the resulting Sobolev spaces
are the ``Sobolev spaces with weights'' considered for instance in
\cite{Kondratiev, KMR1, MazyaRossmann, Rossmann}. Let $f > 0$ be a
smooth function on a domain $\Omega$, we then define the
\emph{$m$th Sobolev space with weight $f$} by
\begin{equation}\label{def.t.c.a0}
    \maK^{m}_{a}(\Omega; f) := \{ u, \
    f^{|\alpha|-a} \pa^{\alpha} u \in L^2(\Omega), \ \
    |\alpha| \le m \}, \quad m \in \NN\cup\{0\},\, a \in \RR.
\end{equation}
Indeed, if $\Omega = \PP \subset \RR^2$ is a polygon, and if we
choose
\begin{equation}
    f(x) = \vartheta(x) = \text{ the {\em distance to
    the non-smooth boundary points of $\PP$}},
\end{equation}
then there is no loss of regularity in the spaces
$\maK^{m}_{a}(\Omega) := \maK^{m}_{a}(\Omega; \vartheta)$
\cite[Theorem~6.6.1]{KMR1}. In this paper, we extend this
regularity result to polyhedral domains in three dimensions,
Theorem~\ref{thm.3D}, with the same choice of the weight (in three
dimensions the weight is the distance to the edges). The analogous
result in arbitrary dimensions leads to topological difficulties
\cite{BNZ-ND, VerVog1}.

Our regularity result requires us first to study the weighted Sobolev
spaces $\maK^{m}_{a}(\Omega) := \maK^{m}_{a}(\Omega; \vartheta)$ where
$\vartheta(x)$ is the distance to the set of singular points on the
boundary. Our approach to Sobolev spaces on polyhedral domains is to
show first that $\maK^{m}_{a}(\Omega)$ is isomorphic to a Sobolev
space on a certain non-compact Riemannian manifold $M$ with smooth
boundary. This non-compact manifold $M$ is obtained from our
polyhedral domain by replacing the Euclidean metric $g_E$ with
\begin{equation}\label{eq.conf.change}
    r_{\PP}^{-2}g_E, \quad r_{\PP} \text{ a smoothing of }
    \vartheta,
\end{equation}
which blows up at the faces of codimension two or higher, that is,
at the set of singular boundary points. (The metric
$r_{\PP}^{-2}g_E$ is Lipschitz equivalent to $\vartheta^{-2}g_E$,
but the latter is not smooth.) The resulting non-compact
Riemannian manifold turns out to be a regular open subset in a
``Lie manifold.'' (see Definition \ref{def.unif.str}, Subsection
\ref{ssec.poly}, and Section \ref{sec.PD} for the precise
definitions). A~{\em Lie manifold} is a compact manifold with
corners $M$ together with a $C^\infty(M)$-module $\VV$ whose
elements are vector fields on $M$. The space $\VV$ must satisfy a
number of axioms, in particular, $\VV$ is required to be closed
under the Lie bracket of vector fields. This property is the
origin of the name \emph{Lie} manifold. The $C^\infty(M)$-module
$\VV$ can be identified with the sections of a vector bundle $A$
over $M$. Choosing a metric on $A$ defines a complete Riemannian
metric on the interior of $M$. See Section~\ref{sec.LI} or
\cite{aln1} for details.

The framework of Lie manifolds is quite convenient for the study
of Sobolev spaces, and in this paper we establish, among other
things, that the main results on the classical Sobolev spaces
remain true in the framework of Lie manifolds. The regular open
sets of Lie manifolds then play in our framework the role played
by smooth, bounded domains in the classical theory.

Let $\PP \subset \RR^n$ be a polyhedral domain. We are
especially interested in describing the spaces
$\maK_{a-1/2}^{m-1/2}(\pa \PP)$ of restrictions to the boundary of
the functions in the weighted Sobolev space $\maK_a^{m}(\PP;
\vartheta) = \maK_a^{m}(\PP; r_{\PP})$ on $\PP$. Using the
conformal change of metric of Equation~\eqref{eq.conf.change}, the
study of restrictions to the boundary of functions in
$\maK_a^{m}(\PP)$ is reduced to the analogous problem on a
suitable regular open subset $\Omega_\PP$ of some Lie manifold.
More precisely, $\maK_a^{m}(\PP) = r_{\PP}^{a - n/2}
H^m(\Omega_\PP)$. A~consequence of this is that
\begin{equation}\label{eq.res}
    \maK_{a-1/2}^{m-1/2}(\pa \PP) = \maK^{m-1/2}_{a-1/2}(\pa \PP;
    \vartheta) = r_{\PP}^{a - n/2}H^{m-1/2}(\pa
    \Omega_\PP).
\end{equation}
(In what follows, we shall usually simply denote
$\maK^{m}_{a}(\PP) := \maK^{m}_{a}(\PP; \vartheta) =
\maK^{m}_{a}(\PP; r_{\PP})$ and $\maK^{m}_{a}(\pa \PP) :=
\maK^{m}_{a}(\pa \PP; \vartheta) = \maK^{m}_{a}(\pa \PP;
r_{\PP})$, where, we recall, $\vartheta(x)$ is the distance from
$x$ to the set of non-smooth boundary points and $r_\PP$ is a
smoothing of $\vartheta$ that satisfies $r_\PP/\vartheta \in [c,
C]$, $c, C > 0$.)

Equation \eqref{eq.res} is one of the motivations to study Sobolev
spaces on Lie manifolds. In addition to the non-compact manifolds
that arise from polyhedral domains, other examples of Lie
manifolds include the Euclidean spaces $\RR^n$, manifolds that are
Euclidean at infinity, conformally compact manifolds, manifolds
with cylindrical and polycylindrical ends, and asymptotically
hyperbolic manifolds. These classes of non-compact
manifolds appear in the study of the Yamabe problem
\cite{lee.parker:87,schoen:84} on compact manifolds, of the Yamabe
problem on asymptotically cylindrical manifolds
\cite{akutagawa.botvinnik}, of analysis on locally symmetric
spaces, and of the positive mass theorem \cite{schoen.yau:79a,
schoen.yau:81, witten:81}, an analogue of the positive mass
theorem on asymptotically hyperbolic manifolds
\cite{andersson.dahl:98}. Lie manifolds also appear in
Mathematical Physics and in Numerical Analysis.
Classes of Sobolev spaces on non-compact manifolds have been
studied in many papers, of which we mention only a few
\cite{ErkipSchrohe, Grubb, Lauter, LauterMoroianu, Maniccia,
Mazzeo, MelroseScattering, meicm, meaps, SchroheSI, SchroheFC,
ScSc, VasyN, VZ} in addition to the works mentioned before. Our
work can also be used to unify some of the various approaches
found in these papers.

Let us now review in more detail the contents of this paper.
A~large part of the technical material in this paper is devoted to
the study of Sobolev spaces on Lie manifolds (with or without
boundary). If $M$ is a {\em compact} manifold with corners, we
shall denote by $\pa M$ the union of all boundary faces of $M$ and
by $M_0 := M \smallsetminus \pa M$ the interior of $M$. We begin
in Section~\ref{sec.LI} with a review of the definition of a
structural Lie algebra of vector fields $\VV$ on a manifold with
corners $M$. This Lie algebra of vector fields will provide the
derivatives appearing in the definition of the Sobolev spaces.
Then we define a Lie manifold as a pair $(M, \VV)$, where $M$ is a
compact manifold with corners and $\VV$ is a structural Lie
algebra of vector fields that is unrestricted in the interior
$M_0$ of $M$. We will explain the above mentioned fact that the
interior of $M$ carries a complete metric $g$. This metric is
unique up to Lipschitz equivalence (or quasi-isometry). We also
introduce in this section Lie manifolds with (true) boundary and,
as an example, we discuss the example of a Lie manifold with true
boundary corresponding to curvilinear polygonal domains. In
Section \ref{sec.SM} we discuss Lie submanifolds, and most
importantly, the global tubular neighborhood theorem. The
proof of this global tubular neighborhood theorem is based on
estimates on the second fundamental form of the boundary, which
are obtained from the properties of the structural Lie algebra of
vector fields. This property distinguishes Lie manifolds from
general manifolds with boundary and bounded geometry, for which a
global tubular neighborhood is part of the definition. In Section
\ref{sec.SS}, we define the Sobolev spaces $W^{s, p}(M_0)$ on the
interior $M_0$ of a Lie manifold $M$, where either $s \in
\NN\cup\{0\}$ and $1 \le p \le \infty$ or $s \in \RR$ and $1 < p <
\infty$. We first define the spaces $W^{s,p}(M_0)$, $s \in
\NN\cup\{0\}$ and $1 \le p \le \infty$, by differentiating with
respect to vector fields in $\VV$. This definition is in the
spirit of the standard definition of Sobolev spaces on $\RR^n$.
Then we prove that there are two alternative, but equivalent ways
to define these Sobolev spaces, either by using a suitable class
of partitions of unity (as in \cite{Shubin, Skrz, Triebel} for
example), or as the domains of the powers of the Laplace operator
(for $p = 2$). We also consider these spaces on open subsets
$\Omega_0 \subset M_0$. The spaces $W^{s, p}(M_0)$, for $s \in
\RR$, $1 < p < \infty$ are defined by interpolation and duality
or, alternatively, using partitions of unity. In Section
\ref{sec.MB}, we discuss regular open subsets $\Omega \subset M$.
In the last two sections, several of the classical results on
Sobolev spaces on smooth domains were extended to the spaces
$W^{s, p}(M_0)$. These results include the density of smooth,
compactly supported functions, the Gagliardo-Nirenberg-Sobolev
inequalities, the extension theorem, the trace theorem, the
characterization of the range of the trace map in the Hilbert
space case ($p = 2$), and the Rellich-Kondrachov compactness
theorem.

In Section \ref{sec.5} we include as an application a regularity
result for strongly elliptic boundary value problems, Theorem
\ref{thm.reg1}. This theorem gives right away the following
result, proved in Section \ref{sec.PD}, which states that there is
no loss of regularity for these problems within weighted Sobolev
spaces.

\begin{theorem}\label{theorem.3Dreg}\
Let $\PP \subset \RR^3$ be a polyhedral domain and $P$ be a strongly
elliptic, second order differential operator with coefficients in
$\CI(\overline{\PP})$. Let $u \in \Kond{1}{a+1}(\PP)$, $u = 0$ on $\pa
\PP$, $a\in\RR$.  If $Pu \in \Kond{m-1}{a-1}(\PP)$, then $u \in
\Kond{m+1}{a+1}(\PP)$ and there exists $C>0$ independent of $u$ such
that
\begin{equation*}
    \|u\|_{\Kond{m+1}{a+1}(\PP)} \le
    C\big(\|Pu\|_{\Kond{m-1}{a-1}(\PP)} +
    \|u\|_{\Kond{0}{a+1}(\PP)}\big), \quad m \in \NN\cup\{0\}.
\end{equation*}
The same result holds for strongly elliptic systems.
\end{theorem}

Note that the above theorem does not constitute a Fredholm (or normal
solvability) result, because the inclusion $\Kond{m+1}{a+1}(\PP) \to
\Kond{0}{a+1}(\PP)$ is {\em not compact}. See also \cite{Kondratiev,
KMR1, MazyaRossmann, Rossmann} and the references therein for similar
results.

In Section \ref{sec.NS}, we obtain a ``non-standard boundary value
problem'' on a smooth domain $\maO$ in weighted Sobolev spaces
with weight given by the distance to the boundary. The boundary
conditions are thus replaced by growth conditions.  Finally, in
the last section, Section \ref{sec.PO}, we obtain mapping
properties for the pseudodifferential calculus
$\Psi_{\VV}^\infty(M)$ defined in \cite{aln2} between our weighted
Sobolev spaces $\rho^s W^{r, p}(M)$. We also obtain a general
elliptic regularity result for elliptic pseudodifferential
operators in $\Psi_{\VV}^\infty(M)$.

\subsubsection*{Acknowledgements:}\ We would like to thank Anna
Mazzucato and Robert Lauter for useful comments.
The first named author wants to thank MSRI,
Berkeley, CA for its hospitality.


\section{Lie manifolds}
\label{sec.LI}

As explained in the Introduction, our approach to the study of
weighted Sobolev spaces on polyhedral domains is based on their
relation to Sobolev spaces on Lie manifolds with true boundary.
Before we recall the definition of a Lie manifold and some of
their basic properties, we shall first look at the following
example, which is one of the main motivations for the theory of
Lie manifolds.

\begin{example}\label{ex.angle}\
Let us take a closer look at the local structure of the Sobolev
space $\maK^m_a(\PP)$ associated to a polygon $\PP$ (recall
\eqref{def.t.c.a0}). Consider $\Omega := \{ (r,\theta)\,|\, 0 <
\theta < \alpha \}$, which models an angle of $\PP$. Then the
distance to the vertex is simply $\vartheta(x) = r$, and the {\em
weighted Sobolev spaces} associated to $\Omega$,
$\maK_a^m(\Omega)$, can alternatively be described as
\begin{equation}\label{def.t.c.a}
    \maK^m_a(\Omega) = \maK^m_a(\Omega; \vartheta):=
    \{ u \in L^2_{loc}(\Omega), \
    r^{-a}(r\pa_r)^i \pa_\theta^ju \in L^2(\Omega), \ \
    i + j \le m \}.
\end{equation}
The point of the definition of the spaces $\maK^m_a(\Omega)$ was
the replacement of the local basis $\{r \pa_x, r \pa_y\}$ with the
local basis $\{r \pa_r, \pa_\theta\}$ that is easier to work with
on the desingularization $\Sigma(\Omega) := [0, \infty) \times [0,
\alpha] \ni (r, \theta)$ of $\Omega$. By further writing $r =
e^t$, the vector field $r\pa_r$ becomes $\pa_t$. Since $dt =
r^{-1}dr$, the space $\maK^m_1(\Omega)$ then identifies with
$H^m(\RR_t\times (0, \alpha))$. The weighted Sobolev space
$\maK^m_1(\Omega)$ has thus become a classical Sobolev space on
the cylinder $\RR\times (0, \alpha)$, as in \cite{Kondratiev}.
\end{example}

The aim of the following definitions is to define such a
desingularisation in general. The desingularisation will carry the
structure of a Lie manifold, defined in the next subsection.

We shall introduce a further, related definition, namely the
definition of a ``Lie submanifolds of a Lie manifold'' in Section
\ref{sec.MB}.

\subsection{Definition of Lie manifolds}\ At first, we want
to recall the definition of manifolds with corners. A
\emph{manifold with corners} is a closed subset $M$ of a
differentiable manifold such that every point $p\in M$ lies in a
coordinate chart whose restriction to $M$ is a diffeomorphism to
$[0,\infty)^k \times \RR^{n-k}$, for some $k = 0, 1, \ldots, n$
depending on $p$. Obviously, this definition includes the property
that the transition map of two different charts are smooth up to
the boundary.  If $k=0$ for all $p \in M$, we shall say that $M$
is a {\em smooth manifold}. If $k\in\{0,1\}$, we shall say that
$M$ is a \emph{smooth manifold with smooth boundary}.

Let $M$ be a {\em compact} manifold with corners. We shall denote
by $\pa M$ the union of all boundary faces of $M$, that is, $\pa
M$ is the union of all points not having a neighborhood
diffeomorphic to $\RR^n$. Furthermore, we shall write $M_0 := M
\smallsetminus \pa M$ for the {\em interior} of~$M$. In order to
avoid confusion, we shall use this notation and terminology only
when $M$ is compact. Note that our definition allows $\pa M$ to be
a smooth manifold, possibly empty.

As we shall see below, a Lie manifold is described by a Lie
algebra of vector fields satisfying certain conditions. We now
discuss some of these conditions.

\begin{definition}\label{def.str.LA}\
A~subspace $\VV\subseteq \Gamma(M; TM)$ of the Lie algebra of all
smooth vector fields on $M$ is said to be a {\em structural Lie
algebra of vector fields on $M$} provided that the following
conditions are satisfied:
\begin{enumerate}[(i)]
\item\ $\VV$ is closed under the Lie bracket of vector fields;
\item\ every $V\in\VV $ is tangent to all boundary hyperfaces of $M$;
\item\ $\CI(M)\VV = \VV$; and
\item\ each point $p\in M$ has a neighborhood $U_p$
such that
\begin{equation*}
    \VV_{U_p} := \{X|_{\overline{U}_p}\,|X\in \VV\} \simeq
    C^\infty(\overline{U}_p)^k
\end{equation*}
in the sense of $C^\infty(\overline{U}_p)$-modules.
\end{enumerate}
\end{definition}

\noindent The condition (iv) in the definition above can be
reformulated as follows:
\begin{enumerate}[(i')]
\setcounter{enumi}3
\item\
For every $p \in M$, there exist a neighborhood $U_p \subset M$ of
$p$ and vector fields $X_1, X_2, \ldots, X_k \in \VV$ with the
property that, for any $Y \in \VV$, there exist functions $f_1,
\ldots, f_k \in \CI(M)$, uniquely determined on $U_p$, such that
\begin{equation}\label{eq.local.basis}
    Y = \sum_{j = 1}^k f_j X_j \quad \text{on} \; U_p.
\end{equation}
\end{enumerate}

We now have defined the preliminaries for the following important
definition.

\begin{definition}\label{def.unif.str}\
A~{\em Lie structure at infinity} on a smooth manifold $M_0$ is a
pair $(M, \VV)$, where $M$ is a compact manifold with interior
$M_0$ and $\VV \subset \Gamma(M; TM)$ is a structural Lie algebra
of vector fields on $M$ with the following property:\ If $p \in
M_0$, then any local basis of $\VV$ in a neighborhood of $p$ is
also a local basis of the tangent space to $M_0$.
\end{definition}

It follows from the above definition that the constant $k$ of
Equation \eqref{eq.local.basis} equals to the dimension $n$ of $M_0$.

A~{\em manifold with a Lie structure at infinity} (or, simply, a
{\em Lie manifold}) is a manifold $M_0$ together with a Lie
structure at infinity $(M, \VV)$ on $M_0$. We shall sometimes
denote a Lie manifold as above by $(M_0, M, \VV)$, or, simply, by
$(M, \VV)$, because $M_0$ is determined as the interior of $M$.
(In \cite{aln1}, only the term ``manifolds with a Lie structure at
infinity'' was used.)

\begin{example}\
If $F \subset TM$ is a sub-bundle of the tangent bundle of a smooth
manifold (so $M$ has no boundary) such that $\VV_F := \Gamma(M; F)$ is
closed under the Lie bracket, then $\VV_F$ is a structural Lie algebra
of vector fields.  Using the Frobenius theorem it is clear that such
vector bundles are exactly the tangent bundles of $k$-dimensional
foliations on $M$, $k=\mathop{\rm rank} F$.  However, $\VV_F$ does not
define a Lie structure at infinity, unless $F=TM$.
\end{example}

\begin{remark}\label{rem.A}\
We observe that Conditions (iii) and (iv) of Definition
\ref{def.str.LA} are equivalent to the condition that $\VV$ be a
projective $\CI(M)$-module. Thus, by the Serre-Swan theorem
\cite{Karoubi}, there exists a vector bundle $A \to M$, unique up
to isomorphism, such that $\VV = \Gamma(M; A)$. Since $\VV$
consists of vector fields, that is $\VV \subset \Gamma(M; TM)$, we
also obtain a natural vector bundle morphism $\varrho_M : A \to
TM$, called the {\em anchor map}. The Condition (ii) of Definition
\ref{def.unif.str} is then equivalent to the fact that $\varrho_M$
is an isomorphism $A\vert_{M_0} \simeq TM_0$ on $M_0$. We will
take this isomorphism to be an identification, and thus we can say
that $A$ is {\em an extension} of $TM_0$ to $M$ (that is, $TM_0
\subset A$).
\end{remark}

\subsection{Riemannian metric\label{ssec.Rm}}
Let $(M_0, M,\VV)$ be a Lie manifold. By definition, a {\em
Riemannian metric on $M_0$ compatible with the Lie structure at
infinity $(M, \VV)$} is a metric $g_0$ on $M_0$ such that, for any
$p \in M$, we can choose the basis $X_1, \dots, X_k$ in
Definition~\ref{def.str.LA} (iv') \eqref{eq.local.basis} to be
orthonormal with respect to this metric everywhere on $U_p\cap
M_0$. (Note that this condition is a restriction only for $p \in
\pa M := M \smallsetminus M_0$.) Alternatively, we will also say
that $(M_0,g_0)$ is a {\em Riemannian Lie manifold}. Any Lie
manifold carries a compatible Riemannian metric, and any two
compatible metrics are bi-Lipschitz to each other.

\begin{remark}\ Using the language of Remark \ref{rem.A}, $g_0$ is a
compatible metric on $M_0$ if, and only if, there exists a metric $g$
on the vector bundle $A \to M$ which restricts to $g_0$ on $TM_0
\subset A$.
\end{remark}

The geometry of a Riemannian manifold $(M_{0},g_{0})$ with a Lie
structure $(M,\VV)$ at infinity has been studied in \cite{aln1}.
For instance, $(M_{0},g_{0})$ is necessarily complete and, if
$\partial M\neq \emptyset$, it is of infinite volume. Moreover,
all the covariant derivatives of the Riemannian curvature tensor
are bounded. Under additional mild assumptions, we also know that
the injectivity radius is bounded from below by a positive
constant, \ie $(M_{0},g_{0})$ is of bounded geometry. (A~{\em
manifold with bounded geometry} is a Riemannian manifold with
positive injectivity radius and with bounded covariant derivatives
of the curvature tensor, see \cite{Shubin} and references
therein).

On a Riemannian Lie manifold $(M_0, M, \VV, g_0)$, the exponential
map $\exp: TM_0 \to M_0$ is well-defined for all $X\in TM_0$ and
extends to a differentiable map $\exp: A \to M$. A~convenient way
to introduce the exponential map is via the geodesic spray, as
done in \cite{aln1}. Similarly, any vector field $X \in \VV =
\Gamma(M; A)$ is integrable and will map any (connected) boundary
face of $M$ to itself. The resulting diffeomorphism of $M_0$ will
be denoted $\psi_X$.

\subsection{Examples}
We include here two examples of Lie manifolds together with compatible
Riemannian metrics. The reader can find more
examples in \cite{aln1,LN1}.

\begin{examples}\label{ex1}\ \\[-5mm]
\begin{enumerate}
\item\ Take $\VV_b$ to be the set of all vector fields tangent to
all faces of a manifold with corners $M$. Then $(M, \VV_b)$ is a
Lie manifold. This generalizes Example~\ref{ex.angle}. See also
Subsection \ref{ssec.poly} and Section \ref{sec.PD}. Let  $r \ge
0$ to be a smooth function on~$M$ that is equal to the distance to
the boundary in a neighborhood of $\pa M$, and is $>0$ outside
$\pa M$ (\ie on~$M_0$). Let $h$ be a smooth metric on $M$, then $g_0
= h + (r^{-1}dr)^2$ is a compatible metric on $M_0$.
\item\ Take $\VV_0$ to be the set of all vector fields vanishing on
all faces of a manifold with corners $M$. Then $(M, \VV_0)$ is a
Lie manifold. If $\pa M$ is a smooth manifold (\ie if $M$ is a
smooth manifold with boundary), then $\VV_0 = r\Gamma(M; TM)$,
where $r$ is as in (a).
\end{enumerate}
\end{examples}

\subsection{$\VV$-differential operators\label{ssec.Diff}}
We are especially interested in the analysis of the differential
operators generated using only derivatives in $\VV$. Let
$\DiffV{*}(M)$ be the algebra of differential operators on $M$
generated by multiplication with functions in $\CI(M)$ and by
differentiation with vector fields $X \in \VV$. The space of order
$m$ differential operators in $\DiffV{*}(M)$ will be denoted
$\DiffV{m}(M )$. A~differential operator in $\DiffV{*}(M)$ will be
called a $\VV$-differential operator.

We can define $\VV$-differential operators acting between sections
of smooth vector bundles $E, F \rightarrow M$, $E, F \subset M
\times \CC^N$ by
\begin{equation}\label{vectorvalued}
    \DiffV{*}(M ;E,F) := e_F M_N(\DiffV{*}(M )) e_E\,,
\end{equation}
where $M_N(\DiffV{*}(M ))$ is the algebra of $N\times N$-matrices
over the ring $\DiffV{*}(M)$, and where
$e_E, e_F \in M_N(\CI(M))$ are the projections onto $E$ and,
respectively, onto $F$. It follows that $\DiffV{*}(M; E):=\DiffV{*}(M; E,E)$
is an algebra. It is also closed under taking
adjoints of operators in $L^2(M_0)$, where the volume form is
defined using a compatible metric $g_0$ on $M_0$.

\subsection{Regular open sets\label{ssec.ROS}}
{\em We assume from now on that $r_{\mathrm{inj}}(M_0)$, the
injectivity radius of $(M_{0}, g_{0})$, is positive.}

One of the main goals of this paper is to prove the results on
weighted Sobolev spaces on polyhedral domains that are needed for
regularity theorems. We shall do that by reducing the study of
weighted Sobolev spaces to the study of Sobolev spaces on ``regular
open subsets'' of Lie manifolds, a class of open sets that plays in
the framework of Lie manifolds the role played by domains with smooth
boundaries in the framework of bounded, open subsets of
$\RR^n$. Regular open subsets are defined below in this subsection.

Let $N \subset M$ be a submanifold of codimension
one of the Lie manifold $(M, \VV)$. Note that this implies that
$N$ is a closed subset of $M$. We shall say that $N$
is a {\em regular} submanifold of $(M,\VV)$ if we can choose a
neighborhood $V$ of $N$ in $M$ and a compatible metric $g_0$ on
$M_0$ that restricts to a product-type metric on $V \cap M_0
\simeq (\pa N_0) \times (-\varepsilon_0, \varepsilon_0)$, $N_0 = N
\smallsetminus \pa N = N \cap M_0$. Such neighborhoods will be
called \emph{tubular neighborhoods}.

In Section \ref{sec.SM}, we shall show that a codimension one
manifold is regular if, and only if, it is a tame submanifold of
$M$; this gives an easy, geometric, necessary and sufficient
condition for the regularity of a codimension one submanifold of
$M$. This is relevant, since the study of manifolds with boundary
and bounded geometry presents some unexpected difficulties
\cite{Schick}.

In the following, it will be important to distinguish properly
between the boundary of a topological subset, denoted by $\patop$,
and the boundary in the sense of manifolds with corners, denoted
simply by $\pa$.

\begin{definition}
Let $(M, \VV)$ be a Lie manifold and $\Omega \subset M$ be an open
subset. We shall say that $\Omega$ is a {\em regular open subset}
in $M$ if, and only if, $\Omega$ is connected, $\Omega$ and
$\overline{\Omega}$ have the same boundary, $\patop \Omega$ (in
the sense of subsets of the topological space $M$), and $\patop
\Omega$ is a regular submanifold of $M$.
\end{definition}

Let $\Omega \subset M$ be a regular open subset. Then
$\overline{\Omega}$ is a compact manifold with corners. The reader
should be aware of the important fact that $\patop\Omega=\patop
\overline \Omega$ is contained in $\pa \overline\Omega$, but in
general $\pa \overline\Omega$ and $\patop\Omega$ are not equal.  The
set $\patop\Omega$ will be called the {\em true boundary} of
$\overline{\Omega}$.  Furthermore, we introduce $\pa_\infty \Omega :=
\pa \overline{\Omega} \cap \pa M$, and call it the {\em boundary at
infinity} of $\overline{\Omega}$.  Obviously, one has $\pa
\overline{\Omega} = \patop \overline{\Omega} \cup \pa_\infty
\overline{\Omega}$.  The true boundary and the boundary at infinity
intersect in a (possibly empty) set of codimension $\geq 2$.  See
Figure~\ref{fig.sigma}.  We will also use the notation $\pa\Omega_0:=
\patop\Omega\cap M_0 = \pa \overline{\Omega}\cap M_0$.

\begin{figure}[htbp]
\begin{center}
\begin{pspicture}(5,0)(11,6)
\newgray{Mgray}{.85}
\newgray{Mogray}{.5}

\rput(5,0){
\pscustom[fillstyle=solid,fillcolor=Mgray,linestyle=none]{
\pscurve(-0.4,0)(1.3,2.7)(3.8,5)(6,5.9) \psline(6,5.9)(6,0)
\psline(5,0)(-0.4,0)}
\pscustom[fillstyle=solid,fillcolor=Mogray,linestyle=dashed]{
\psecurve[linestyle=solid](-0.4,0)(1.3,2.7)(3.8,5)(6,5.9)
\pscurve[linestyle=dashed](3.8,5)(4.8,3)(3.2,1.8)(1.3,2.7) }
\pscurve[linewidth=1.3pt](-0.4,0)(1.3,2.7)(3.8,5)(6,5.9)
\rput(3.2,3.3){\psframebox[linecolor=white,fillstyle=solid,fillcolor=white]
{$\Omega$}}
\rput(5.2,2.4){\psframebox[linecolor=white,fillstyle=solid,fillcolor=white]
{$\patop\overline{\Omega}$}}
\rput(2.4,4.8){\psframebox[linecolor=white,fillstyle=solid,fillcolor=white]
{$\partial_\infty\Omega$}}
\rput(1.6,1){\psframebox[linecolor=white,fillstyle=solid,fillcolor=white]
{$M_0$}} }
\end{pspicture}
\end{center}
\caption{A~regular open set $\Omega$. Note that the interior of
$\pa_\infty \overline{\Omega}$ is contained in $\Omega$, but the
true boundary $\patop \Omega =\patop \overline\Omega$ is not
contained in $\Omega$} \label{fig.sigma}
\end{figure}

The space of restrictions to $\Omega$ or $\overline{\Omega}$ of
order $m$ differential operators in $\DiffV{*}(M)$ will be denoted
$\DiffV{m}(\Omega)$, respectively $\DiffV{m}(\overline{\Omega})$.
Similarly, we shall denote by $\VV(\Omega)$ the space of
restrictions to $\overline{\Omega}$ of vector fields in $\VV$, the
structural Lie algebra of vector fields on $M$.

Let $F\subset \pa \Omega$ be any boundary hyperface of $\overline
\Omega$ of codimension $1$.  Such a face is either contained in
$\patop \overline{\Omega}$ or in $\pa_\infty \overline{\Omega}$.
If $F\subset \pa_\infty \Omega$, then the restrictions of all
vector fields in $\VV$ to $F$ are tangent to $F$. However, if
$F\subset \patop \overline\Omega$ the regularity of the boundary
implies that there are vector fields in $\VV$ whose restriction to
$F$ is not tangent to $F$. In particular, the true boundary
$\patop \overline\Omega$ of $\overline{\Omega}$ is uniquely
determined by $(\overline{\Omega}, \VV(\Omega))$, and hence so is
$\Omega = \overline{\Omega} \smallsetminus \patop
\overline\Omega$. We therefore obtain a one-to-one correspondence
between Lie manifolds with true boundary and regular open subsets
(of some Lie manifold $M$).

Assume we are given $\Omega$, $\overline{\Omega}$ (the closure in
$M$), and $\VV(\Omega)$, with $\Omega$ a regular open subset of
some Lie manifold $(M, \VV)$. In the cases of interest, for
example if $\patop \overline\Omega$ is a tame submanifold of $M$
(see Subsection \ref{ssec.tame} for the definition of tame
submanifolds), we can replace the Lie manifold $(M, \VV)$ in which
$\Omega$ is a regular open set with a Lie manifold $(N, \WW)$
canonically associated to $(\Omega, \overline{\Omega},
\VV(\Omega))$ as follows.  Let $N$ be obtained by gluing two
copies of $\overline{\Omega}$ along $\patop \overline\Omega$, the
so-called {\em double} of $\overline{\Omega}$, also denoted
$\overline{\Omega}^{db} = N$. A~smooth vector field $X$ on
$\overline{\Omega}^{db}$ will be in $\WW$, the structural Lie
algebra of vector fields $\WW$ on $\overline{\Omega}^{db}$ if, and
only if, its restriction to each copy of $\overline{\Omega}$ is in
$\VV(\Omega)$. Then $\Omega$ will be a regular open set of the Lie
manifold $(N, \WW)$. For this reason, the pair
$(\overline{\Omega}, \VV(\Omega))$ will be called a {\em Lie
manifold with true boundary}. In particular, the true boundary of
a Lie manifold with true boundary is a tame submanifold of the
double. The fact that the double is a Lie manifold is justified in
Remark \ref{rem.tub}.

\subsection{Curvilinear polygonal domains\label{ssec.poly}}
We conclude this section with a discussion of a curvilinear
polygonal domain $\PP$, an example that generalizes Example
\ref{ex.angle} and is one of the main motivations for considering
Lie manifolds. To study function spaces on $\PP$, we shall
introduce a ``desingularization'' $(\Sigma(\PP), \kappa)$ of $\PP$
(or, rather, of~$\overline{\PP}$), where $\Sigma(\PP)$ is a
compact manifold with corners and $\kappa : \Sigma(\PP) \to
\overline{\PP}$ is a continuous map that is a diffeomorphism from
the interior of $\Sigma(\PP)$ to $\PP$ and maps the boundary of
$\Sigma(\PP)$ onto the boundary of $\PP$.

Let us denote by $B^k$ the open unit ball in $\RR^k$.

\begin{definition}\label{def.curv.poly}\
An open, connected subset $\PP \subset M$ of a two dimensional
manifold~$M$ will be called a {\em curvilinear polygonal domain}
if, by definition, $\overline{\PP}$ is compact and for every point
$p \in \pa \PP$ there exists a diffeomorphism $\phi_p : V_p \to
B^2$, $\phi_p(p) = 0$, defined on a neighborhood $V_p \subset M$
such that
\begin{equation}\label{eq.phi.2D}
    \phi_j(V_p \cap \PP) = \{(r \cos \theta, r \sin \theta),
    \, 0 < r < 1,\, 0 < \theta < \alpha_p\}\,,
    \quad \alpha_p \in (0, 2\pi).
\end{equation}
\end{definition}

A~point $p \in \pa \PP$ for which $\alpha_p \neq \pi$ will be
called a {\em vertex of $\PP$}. The other points of $\pa \PP$ will
be called {\em smooth boundary points}. It follows that every
curvilinear polygonal domain has finitely many vertices and its
boundary consists of a finite union of smooth curves $\gamma_j$
(called the {\em edges} of $\PP$) which have no other common
points except the vertices. Moreover, every vertex belongs to
exactly two edges.

Let $\{P_1, P_2, \ldots, P_k\} \subset \overline{\PP}$ be the
vertices of $\PP$. The cases $k = 0$ and $k = 1$ are also allowed.
Let $V_{j} := V_{P_j}$ and $\phi_j := \phi_{P_j} : V_j \to B^2$ be
the diffeomorphisms defined by Equation \eqref{eq.phi.2D}. Let
$(r, \theta) : \RR^2 \smallsetminus \{(0,0)\} \to (0, \infty)
\times [0, 2\pi)$ be the polar coordinates. We can assume that the
sets $V_j$ are disjoint and define $r_j(x) = r(\phi_j(x))$ and
$\theta_j(x) = \theta(\phi_j(x))$.

The {\em desingularization} $\Sigma(\PP)$ of $\PP$ will replace each
of the vertices $P_j$, $j = 1, \ldots, k$ of $\PP$ with a segment of
length $\alpha_j = \alpha_{P_j} > 0$. Assume that $\PP \subset
\RR^2$. We can realize $\Sigma(\PP)$ in $\RR^3$ as follows. Let
$\psi_j$ be smooth functions supported on $V_j$ with $\psi_j =1$ in a
neighborhood of $P_j$.
\begin{equation*}
    \Phi : \overline{\PP} \smallsetminus \{P_1, P_2, \ldots, P_k\} \to
    \RR^2 \times \RR, \quad \Phi(p) = \big(\,p\,,\, \sum_{j} \psi_j(p)
    \theta_j(p)\,\big).
\end{equation*}
Then $\Sigma(\PP)$ is (up to a diffeomorphism) the closure of
$\Phi(\PP)$ in $\RR^3$. The desingularization map is $\kappa(p,
z) = p$.

The structural Lie algebra of vector fields $\VV(\PP)$ on
$\Sigma(\PP)$ is given by (the lifts of) the smooth vector fields
$X$ on $\overline{\PP} \smallsetminus \{P_1, P_2, \ldots, P_k\}$
that, on $V_j$, can be written as
\begin{equation*}
    X = a_{r}(r_j, \theta_j) r_j \pa_{r_j} + a_{\theta} (r_j,
    \theta_j) \pa_{\theta_j},
\end{equation*}
with $a_r$ and $a_\theta$ smooth functions of $(r_j, \theta_j)$,
$r_j \ge 0$. Then $(\Sigma(\PP), \VV(\PP))$ is a
Lie manifold with true boundary.

To define the structural Lie algebra of vector fields on
$\Sigma(\PP)$, we now choose a smooth function $r_{\PP} : \PP \to
[0, \infty)$ with the following properties
\begin{enumerate}[(i)]
\item\ $r_{\PP}$ is continuous on $\overline{\PP}$,
\item\ $r_{\PP}$ is smooth on $\PP,$
\item\ $r_{\PP}(x) > 0$ on $\overline{\PP} \smallsetminus
\{P_1, P_2, \ldots, P_k\}$,
\item\ $r_{\PP}(x) = r_j(x)$ if $x \in V_j$.
\end{enumerate}

Note that $r_{\PP}$ lifts to a {\em smooth positive function} on
$\Sigma(\PP)$. Of course, $r_{\PP}$ is determined only up to a
smooth positive function $\psi$ on $\Sigma(\PP)$ that equals to $1$ in
a neighborhood of the vertices.

\begin{definition}\label{def.can.weight}
A~function of the form $\psi r_\PP$, with $\psi \in
\CI(\Sigma(\PP))$, $\psi>0$ will be called a {\em canonical weight
function of $\PP$}.
\end{definition}

In what follows, we can replace $r_{\PP}$ with any canonical
weight function. Canonical weight functions will play an important
role again in Section \ref{sec.PD}. Canonical weights are
example of ``admissible weights,'' which will be used to define
weighted Sobolev spaces.

Then an alternative definition of $\VV(\PP)$ is
\begin{equation}\label{eq.vf1}
    \VV(\PP) := \{\, r_{\PP}\,(\psi_1 \pa_1 + \psi_2 \pa_2)\,\}, \quad
    \psi_1, \psi_2 \in \CI(\Sigma(\PP)).
\end{equation}
Here $\pa_1$ denotes the vector field corresponding to the
derivative with respect to the first component. The vector field
$\pa_2$ is defined analogously. In particular,
\begin{equation}\label{eq.vf2}
    r_{\PP} (\pa_j r_{\PP}) = r_{\PP}
    \frac{\pa r_{\PP}}{\pa x_j} \in \CI(\Sigma(\PP)),
\end{equation}
which is useful in establishing that $\VV(\PP)$ is a Lie algebra.
Also, let us notice that both $\{r_{\PP}\pa_1,r_{\PP}\pa_2\} $ and
$\{r_{\PP}\pa_{{r_\PP}},\pa_\theta\}$ are local bases for $\VV(\PP)$ on
$V_{j}$. The transition functions lift to smooth functions on $\Sigma(\PP)$
defined in a neighborhood of $\kappa^{-1}(P_j)$, but cannot be
extended to smooth functions defined in a neighborhood of $P_j$ in
$\overline{\PP}$.

Then $\patop \Sigma(\PP)$, the true boundary of $\Sigma(\PP)$,
consists of the disjoint union of the edges of $\PP$ (note that
the interiors of these edges have disjoint closures in
$\Sigma(\PP)$). Anticipating the definition of a Lie submanifold
in Section \ref{sec.SM}, let us notice that $\patop \Sigma(\PP)$
is a Lie submanifold, where the Lie structure consists of the
vector fields on the edges that vanish at the end points of the
edges.

The function $\vartheta$ used to define the Sobolev spaces
$\maK_a^m(\PP) := \maK_a^m(\PP; \vartheta)$ in
Equation~\eqref{def.t.c.a0} is closely related to the function
$r_\PP$. Indeed, $\vartheta(x)$ is the distance from $x$ to the
vertices of $\PP$. Therefore $\vartheta/r_\PP$ will extend to a
continuous, nowhere vanishing function on $\Sigma(\PP)$, which
shows that
\begin{equation}\label{eq.same.Sobo}
    \maK_a^m(\PP; \vartheta) = \maK_a^m(\PP; r_{\PP}).
\end{equation}

If $P$ is an order $m$ differential operator with smooth
coefficients on $\RR^2$ and $\PP \subset \RR^2$ is a polygonal
domain, then $r_{\PP}^m P \in \DiffV{m}(\Sigma(\PP))$, by Equation
\eqref{eq.vf1}. However, in general, $r_{\PP}^m P$ will not define
a smooth differential operator on $\overline{\PP}$.

\section{Submanifolds\label{sec.SM}}

In this section we introduce various classes of submanifolds of a
Lie manifold. Some of these classes were already mentioned in the
previous sections.

\subsection{General submanifolds}
We first introduce the most general class of  submanifolds of a
Lie manifold.

We first fix some notation. Let $(M_0, M, \VV)$ and $(N_0,
N, \WW)$ be Lie manifolds. We know that there exist vector
bundles $A \to M$ and $B \to N$ such that $\VV \simeq \Gamma(M;
A)$ and $\WW \simeq \Gamma(N; B)$, see Remark \ref{rem.A}. We can
assume that $\VV = \Gamma(M; A)$ and $\WW = \Gamma(N; B)$
and write $(M, A)$ and $(N, B)$ instead of $(M_0, M,
\VV)$ and  $(N_0, N, \WW)$.

\begin{definition}\label{def.submanif}\
Let $(M, A)$ be a Lie manifold with anchor map $\varrho_M : A \to
TM$. A~Lie manifold $(N, B)$ is called a {\em Lie submanifold} of
$(M, A)$ if
\begin{enumerate}[{\rm (i)}]
\item\ $N$ is a closed submanifold of $M$ (possibly with corners, no
transversality at the boundary required),
\item\ $\pa N = N \cap \pa M$ (that is, $N_0 \subset M_0$,
$\partial N \subset \partial M$), and
\item\ $B$ is a sub vector bundle of $A|_N$, and
\item\label{c.BA.rel}\ the restriction
of $\varrho_M$ to $B$ is the anchor map of $B \to N$.
\end{enumerate}
\end{definition}

\begin{remark}\ An alternative form of Condition
\eqref{c.BA.rel} of the above definition is
\begin{multline}\label{eq.alt.f}
    \WW = \Gamma(N; B) = \{X\vert_N\,|\,  X \in \Gamma(M; A)
    \mbox{ and } X \vert_N \text{ tangent to } N\} \\
    = \{X\in \Gamma(N; A|_N)\,|\,
    \varrho_M\circ X \in \Gamma(N; TN)\}.
\end{multline}
\end{remark}

We have the following simple corollary that justifies Condition
\eqref{c.BA.rel} of Definition~\ref{def.submanif}.

\begin{corollary}\label{cor.res.met}\
Let $g_0$ be a metric on $M_0$ compatible with the Lie structure
at infinity on $M_0$. Then the restriction of $g_0$ to $N_0$ is
compatible with the Lie structure at infinity on $N_0$.
\end{corollary}

\begin{proof}\ Let $g$ be a metric on $A$ whose restriction to
$TM_0$ defines the metric $g_0$. Then $g$ restricts to a metric $h$
on $B$, which in turn defines a metric $h_0$ on $N_0$. By
definition, $h_0$ is the restriction of $g_0$ to $N_0$.
\end{proof}

We thus see that any submanifold (in the sense of the above
definition) of a Riemannian Lie manifold is itself a Riemannian
Lie manifold.

\subsection{Second fundamental form}
We define the \emph{$A$-normal bundle} of the Lie submanifold $(N,
B)$ of the Lie manifold $(M, A)$ as $\smb = (A|_N)/B$ which is a
bundle over $N$. Then the anchor map $\varrho_M$ defines a
map $\nu^A \to (TM|_N)/TN$, called the \emph{anchor map of
$\nu^A$}, which is an isomorphism over $N_0$.

We denote the Levi-Civita-connection on $A$ by $\nabla^A$ and the
Levi-Civita connection on $B$ by $\nabla^B$ \cite{aln1}. Let
$X,Y,Z\in \WW = \Gamma(N; B)$ and $\tilde X, \tilde Y,\tilde Z \in
\VV = \Gamma(M; A)$ be such that $X = \tilde X\vert_N$, $Y
=\tilde Y\vert_N$, $Z =\tilde Z\vert_N$. Then $\nabla^A_{\tilde
X} \tilde Y|_N$ depends only on $X, Y \in \WW = \Gamma(N; B)$ and
will be denoted $\nabla^A _XY$ in what follows.
Furthermore, the Koszul formula gives
\begin{align*}
  2 g(\tilde Z,\nabla^A_{\tilde Y} {\tilde X})=&
  \partial_{\varrho_M({\tilde X})}g({\tilde Y},{\tilde Z})
  +\partial_{\varrho_M({\tilde Y})}g({\tilde Z},{\tilde X})
  -\partial_{\varrho_M({\tilde Z})}g({\tilde X},{\tilde Y})\\
  &{}-g([{\tilde X},{\tilde Z}],{\tilde Y})-g([{\tilde Y},{\tilde
  Z}], {\tilde X})-g([{\tilde X},{\tilde Y}],{\tilde Z}),\\
  2 g(Z,\nabla^B_Y X)&=
  \partial_{\varrho_M(X)}g(Y,Z)
  +\partial_{\varrho_M(Y)}g(Z,X)
  -\partial_{\varrho_M(Z)}g(X,Y)\\
  &{}-g([X,Z],Y)-g([Y,Z],X)-g([X,Y],Z).
\end{align*}
As this holds for arbitrary sections $Z$ of $\Gamma(N; B)$ with
extensions $\tilde Z$ on $\Gamma(M;A)$, we see that $\nabla^B_X Y$
is the tangential part of $\nabla^A_X Y|_N$.

The normal part of $\nabla^A$ then gives rise to the \emph{second
fundamental form} $\II$ defined as
\begin{equation*}
    \II : \WW \times \WW \to \Gamma(\smb), \quad \II(X,Y) :=
    \nabla^A_X Y - \nabla^B_X Y.
\end{equation*}
The Levi-Civita connections $\nabla^A$ and $\nabla^B$ are torsion
free, and hence $\II$ is symmetric because
\begin{equation*}
    \II(X,Y)-\II(Y,X)=[\tilde X,\tilde Y]\vert_N-[X,Y] = 0.
\end{equation*}

A~direct computation reveals also that $\II(X,Y)$ is tensorial in
$X$, and hence, because of the symmetry, it is also tensorial in
$Y$. (``Tensorial'' here means $\II(fX, Y) = f\II(X, Y) = \II(X,
fY)$, as usual.) Therefore the second fundamental form is a
vector bundle morphism $\II : B \otimes B \to \smb$, and the
endomorphism at $p\in M$ is denoted by $\II_p:B_p\otimes B_p\to
A_p$. It then follows from the compactness of $N$ that
\begin{equation*}
    \|\II_p(X_p, Y_p)\|\leq C \|X_p\|\,\|Y_p\|,
\end{equation*}
with a constant $C$ independent of $p \in N$. Clearly, on the
interior $N_0\subset M_0$ the second fundamental form coincides
with the classical second fundamental form.

\begin{corollary}\label{cor.II.bdd}
Let $(N,B)$ be a submanifold of $(M,A)$ with a compatible
metric. Then the (classical) second fundamental form of $N_0$ in
$M_0$ is uniformly bounded.
\end{corollary}

\subsection{Tame submanifolds\label{ssec.tame}}
We now introduce tame manifolds. Our main interest in tame
manifolds is the global tubular neighborhood theorem, Theorem
\ref{thm.tub.neigh}, which asserts that a tame submanifold of a
Lie manifold has a tubular neighborhood in a strong sense. In
particular, we will obtain that a tame submanifold of codimension
one is regular. This is interesting because being tame is an
algebraic condition that can be easily verified by looking at the
structural Lie algebras of vector fields. On the other hand, being
a regular submanifold is an analytic condition on the metric that
may be difficult to check directly.

\begin{definition}\label{def.submanif.ii}\
Let $(N,B)$ be a Lie submanifold of the  Lie manifold
$(M,A)$ with anchor map $\varrho_M : A \to TM$. Then $(N, B)$ is
called a {\em tame submanifold of $M$} if $T_p N$ and
$\varrho_M(A_p)$ span $T_p M$ for all $p \in \partial N$.
\end{definition}

Let $(N, B)$ be a tame submanifold of the Lie manifold $(M, A)$.  Then
the anchor map $\varrho_M : A \to TM$ defines an isomorphism from
$A_p/B_p$ to $T_pM/T_pN$ for any $p\in N$. In particular, the anchor
map $\varrho_M$ maps $B^\perp$, the orthogonal complement of $B$ in
$A$, injectively into $\varrho_M(A)\subset TM$. For any boundary face
$F$ and $p\in F$ we have $\varrho_M(A_p)\subset T_pF$.  Hence, for any
$p\in N\cap F$, the space $T_pM$ is spanned by $T_pN$ and $T_pF$.  As
a consequence, $N \cap F$ is a submanifold of $F$ of codimension $\dim
M -\dim N$. The codimension of $N \cap F$ in $F$ is therefore
independent of $F$, in particular independent of the dimension of $F$.

\begin{examples}\ \\[-5mm]
\begin{enumerate}[{\rm (1)}]
\item Let $M$ be any compact manifold (without boundary).
Fix a $p\in M$. Let $(N,B)$ be a manifold with a Lie structure at
infinity. Then $(N_0\times \{p\},N\times \{p\},B)$ is a tame
submanifold of $(N_0\times M,N\times M,B\times TM)$.
\item If $\partial N\neq \emptyset$, the diagonal
$N$ is a submanifold of $N\times N$, but not a tame submanifold.
\item Let $N$ be a submanifold with corners of $M$ such that $N$ is
transverse to all faces of $M$. We endow these manifolds with the
$b$-structure at infinity $\VV_b$ (see Example \ref{ex1} (i)).
Then $(N, \VV_b)$ is a tame Lie submanifold of $(M, \VV_b)$.
\item A~regular submanifold (see section~1)
is a also a tame submanifold.
\end{enumerate}
\end{examples}

We now prove the main theorem of this section. Note that this
theorem is not true for a general manifold of bounded geometry
with boundary (for a manifold with bounded geometry and boundary,
the existence of a global tubular neighborhood of the boundary is
part of the definition, see \cite{Schick}).

\begin{theorem}[Global tubular neighborhood theorem]
\label{thm.tub.neigh}  Let $(N,B)$ be a tame submanifold of the
Lie manifold $(M,A)$. For $\epsilon>0$, let $(\nu^A)_\epsilon$ be
the set of all vectors normal to $N$ of length smaller than
$\epsilon$. If $\epsilon > 0$ is sufficiently small, then the
normal exponential map $\exp^\nu$ defines a diffeomorphism from
$(\nu^A)_\epsilon$ to an open neighborhood $V_\epsilon$ of~$N$ in
$M$. Moreover, $dist(\exp^\nu(X),N)=|X|$ for $|X|<\ep$.
\end{theorem}

\begin{proof}\
Recall from \cite{aln1} that the exponential map $\exp: TM_0 \to
M_0$ extends to a map $\exp: A \to M$. The definition of the
normal exponential function $\exp^\nu$ is obtained by identifying
the quotient bundle $\nu^A$ with $B^{\perp}$, as discussed
earlier. This gives
\begin{equation*}
    \exp^\nu : (\nu^A)_\epsilon \to M.
\end{equation*}
The differential $d\exp^\nu$ at $0_p\in \nu^A_p$, $p\in N$ is the
restriction of the anchor map to $B^\perp\cong \nu^A$,
hence any point $p\in N$ has a neighborhood $U(p)$ and
$\tau_p>0$ such that
\begin{equation}\label{map.res}
    \exp^\nu:(\nu^A)_{\tau_p}|_{U_p}\to M
\end{equation}
is a diffeomorphism onto its image. By compactness $\tau_p\geq
\tau>0$. Hence, $\exp^\nu$ is a local diffeomorphism of
$(\nu^A)_\tau$ to a neighborhood of~$N$ in $M$. It remains to show
that it is injective for small $\epsilon\in(0,\tau)$.

Let us assume now that there is no $\epsilon>0$ such that the
theorem holds. Then there are sequences $X_i,Y_i\in \nu^A$, $i
\in \NN$, $X_i\neq Y_i$ such that $\exp^\nu X_i=\exp^\nu Y_i$ with
$|X_i|,|Y_i|\to 0$ for $i \to \infty$. After taking a subsequence
we can assume that the basepoints $p_i$ of $X_i$ converge to
$p_\infty$ and the basepoints $q_i$ of $Y_i$ converge to
$q_\infty$. As the distance in $M$ of $p_i$ and $q_i$ converges to
$0$, we conclude that $p_\infty=q_\infty$. However, $\exp^\nu$ is
a diffeomorphism from $(\nu^A)_{\tau}|_{U(p_\infty)}$ into a
neighborhood of $U(p_\infty)$. Hence, we see that $X_i=Y_i$ for
large $i$, which  contradicts the assumptions.
\end{proof}

We now prove that every tame codimension one Lie submanifold is
regular.

\begin{proposition}\label{prop.prod.struct}\
Let $(N,B)$ be a tame submanifold of codimension one of $(M,A)$.
We fix a unit length section $X$ of $\nu^A$.
Theorem~\ref{thm.tub.neigh} states that
\begin{eqnarray*}
    \exp^\nu:(\nu^A)_\epsilon\cong N\times (-\epsilon,\epsilon)&\to&
    \{x\,|\,d(x,N)<\epsilon\} =: V_\epsilon \cr
     (p,t)& \mapsto & \exp\bigl(tX(p)\bigr)
\end{eqnarray*}
is a diffeomorphism for small $\ep>0$.
Then $M_0$ carries a compatible
metric $g_0$ such that $(\exp^\nu)^*g_0$ is a product metric, \ie
$(\exp^\nu)^*g_0 = g_N + dt^2$ on $N \times
(-\epsilon/2,\epsilon/2)$.
\end{proposition}

\begin{proof}\
Choose any compatible metric $g_1$ on $M_0$. Let $g_2$ be a metric
on $U_\epsilon$ such that $(\exp^\nu)^*g_2 = g_1|_N + dt^2$ on
$N\times (-\epsilon, \epsilon)$. Let $d(x) := dist(x, N)$. Then
\begin{equation*}
    g_0 =(\chi\circ d)\; g_1 +(1-\chi\circ d)\; g_2,
\end{equation*}
has the desired properties, where the cut-off function $\chi :
\RR \to [0,1]$ is $1$ on $(-\epsilon/2, \epsilon/2)$ and has support
in $(-\epsilon,\epsilon)$, and satisfies $\chi(-t)=\chi(t)$.
\end{proof}

The above definition shows that any tame submanifold of
codimension~$1$ is a regular submanifold. Hence, the concept of a
tame submanifold of codimension~$1$ is the same as that of a
regular submanifolds. We hence obtain a new criterion for deciding
that a given domain in a Lie manifold is regular.

\begin{proposition}
Assume the same conditions as the previous proposition.  Then
$d\exp^\nu\left(\frac{\pa}{\pa t}\right)$ defines a smooth vector
field on $V_{\ep/2}$. This vector field can be extended smoothly to a
vector field $Y$ in $\VV$. The restriction of $A$ to $V_{\ep/2}$
splits in the sense of smooth vector bundles as $A=A_1\oplus A_2$
where $A_1|_N=\nu^A$ and $A_2|_N=B$. This splitting is parallel in the
direction of $Y$ with respect to the Levi-Civita connection of the
product metric $g_0$, i.e.\ if $Z$ is a section of $A_i$, then
$\nabla_{Y}Z$ is a section of $A_i$ as well.
\end{proposition}

\begin{proof}
Because of the injectivity of the normal exponential map, the
vector field $Y_1:=d\exp^\nu\left(\frac{\pa}{\pa t}\right)$ is
well-defined, and the diffeomorphism property implies smoothness
on $V_\ep$.  At first, we want to argue that $Y_1\in \VV(V_\ep)$.
Let $\pi:S(A)\to M$ be the bundle of unit length vectors in $A$.
Recall from \cite{aln1}, section 1.2 that $S(A)$ is naturally a
Lie manifold, whose Lie structure is given by the \emph{thick
pullback} $\pi^\#(A)$ of $A$. Now the flow lines of $Y_1$ are
geodesics, which yield in coordinates solutions to a second order
ODE in $t$.  In \cite{aln1}, section~3.4 this ODE was studied on
Lie manifolds.  The solutions are integral lines of the geodesic
spray $\sigma:S(A)\to f^\#(A)$. As the integral lines of this flow
stay in $S(A)\subset A$ and as they depend smoothly on the initial
data and on $t$, we see that $Y_1$ is a smooth section of constant
length one of~$A|_{V_\ep}$.

Multiplying with a suitable cutoff-function with support in
$V_\ep$ one sees that we obtain the desired extension $Y\in \VV$.
Using parallel transport in the direction of $Y$, the splitting
$A|_N=\nu^A\oplus TN$ extends to a small neighborhood of $N$. This
splitting is clearly parallel in the direction of $Y$.
\end{proof}

\begin{remark}\label{rem.tub}
Let $N \subset M$ be a tame submanifold of the Lie manifold $(M,
\VV)$ and $Y \in \VV$ as above. If $Y$ has length one in a
neighborhood of $N$ and is orthogonal to $N$, then $V : =
\bigcup_{|t| < \epsilon} \phi_t(N)$ will be a tubular neighborhood
of $N$. According to the previous proposition the restriction of
$A \to M$ to $V$ has a natural product type decomposition. This
justifies, in particular, that the double of a Lie manifold with
boundary is again a Lie manifold, and that the Lie structure
defined on the double satisfies the natural compatibility
conditions with the Lie structure on a Lie manifold with boundary.
\end{remark}

\section{Sobolev spaces\label{sec.SS}}

In this section we study Sobolev spaces on Lie manifolds without
boundary. These results will then be used to study Sobolev spaces
on Lie manifolds with true boundary, which in turn, will be used
to study weighted Sobolev spaces on polyhedral domains. The goal
is to extend to these classes of Sobolev spaces the main results
on Sobolev spaces on smooth domains.
\smallskip

{\bf Conventions.}\ {\em Throughout the rest of this paper, $(M_0,
M, \VV)$ will be a fixed Lie manifold. We also fix a compatible
metric $g$ on $M_0$, \ie a metric compatible with the Lie
structure at infinity on $M_0$, see Subsection \ref{ssec.Rm}. To
simplify notation we denote the compatible metric by $g$ instead
of the previously used $g_0$. By $\Omega$ we shall denote an open
subset of $M$ and $\Omega_0=\Omega\cap M_0$. The letters $C$ and
$c$ will be used to denote possibly different constants that may
depend only on $(M_0,g)$ and its Lie structure at infinity
$(M,\VV)$.}
\smallskip

We shall denote the volume form (or measure) on $M_0$ associated
to $g$ by $d\vol_g(x)$ or simply by $dx$, when there is no danger
of confusion. Also, we shall denote by $L^p(\Omega_0)$ the
resulting $L^p$-space on $\Omega_0$ (\ie defined with respect to
the volume form $dx$). These spaces are independent of the choice
of the compatible metric $g$ on $M_0$, but their norms, denoted by
$\|\; \cdot \;\|_{L^p}$, do depend upon this choice, although this
is not reflected in the notation. Also, we shall use the fixed
metric $g$ on $M_0$ to trivialize all density bundles. Then the
space $\maD'(\Omega_0)$ of distributions on $\Omega_0$ is defined,
as usual, as the dual of $\CIc(\Omega_0)$. The spaces
$L^p(\Omega_0)$ identify with spaces of distributions on
$\Omega_0$ via the pairing
\begin{equation*}
    \langle u, \phi \rangle = \int_{\Omega_0} u(x) \phi(x)dx,\,
    \quad \,\text{where }\, \phi \in \CIc(\Omega_0)
    \,\text{ and }\, u \in L^p(\Omega_0).
\end{equation*}

\subsection{Definition of Sobolev spaces using vector fields and
connections} We shall define the Sobolev spaces
$W^{s,p}(\Omega_0)$ in the following two cases:
\begin{itemize}
\item\ $s\in \NN\cup\{0\}$, $1\leq p\leq \infty$, and arbitrary
open
sets $\Omega_0$ or
\item\ $s\in \RR$, $1<p<\infty$, and $\Omega_0=M_0$.
\end{itemize}
We shall denote $W^{s, p}(\Omega) = W^{s, p}(\Omega_0)$ and $W^{s,
p}(M) = W^{s, p}(M_0)$. If $\Omega$ is a regular open set, then
$W^{s, p}(\overline{\Omega}) = W^{s, p}(\Omega_0)$. In the case
$p=2$, we shall often write $H^s$ instead of $W^{s, 2}$.  We shall
give several definitions for the spaces $W^{s, p}(\Omega_0)$ and
show their equivalence. This will be crucial in establishing the
equivalence of various definitions of weighted Sobolev spaces on
polyhedral domains. The first definition is in terms of the
Levi-Civita connection $\nabla$ on $TM_0$. We shall denote also by
$\nabla$ the induced connections on tensors (\ie on tensor
products of $TM_0$ and $T^*M_0$).

\begin{definition}[$\nabla$-definition of Sobolev
spaces]\label{def.Sobolev}\ The Sobolev space $W^{k,
p}(\Omega_0)$, $k \in \NN\cup\{0\}$, is defined as the space of
distributions $u$ on $\Omega_0 \subset M_0$ such that
\begin{equation}\label{eq.norm.Wkp'}
    \|u\|_{\nabla, \Wkp}^p : = \sum_{l = 1}^k \int_{\Omega_0}
    |\nabla^l u(x)|^pdx < \infty \,, \quad 1 \le p < \infty.
\end{equation}
For $p = \infty$ we change this definition in the obvious way,
namely we require that,
\begin{equation}\label{eq.norm.Wkinf'}
    \|u\|_{\nabla, \Wkinf} : = \sup |\nabla^l u(x)|
    < \infty \,, \quad 0 \le l \le k.
\end{equation}
\end{definition}

We introduce an alternative definition of Sobolev spaces.

\begin{definition}[vector fields definition of Sobolev spaces]
Let again $k\in \NN\cup\{0\}$. Choose a finite set of vector
fields $\maX$ such that $\CI(M) \maX = \VV$. This condition is
equivalent to the fact that the set $\{X(p), X \in \maX\}$
generates $A_p$ linearly, for any $p \in M$. Then the system
$\maX$ provides us with the norm
\begin{equation}\label{eq.norm.Wkp}
    \|u\|_{\maX, \Wkp}^p : = \sum \|X_{1} X_{2} \ldots X_{l} u
    \|_{L^p}^p\,,\quad 1 \le p < \infty,
\end{equation}
the sum being over all possible choices of $0 \le l \le k$ and all
possible choices of not necessarily distinct vector fields $X_1,
X_2, \ldots, X_l \in \maX$. For $p = \infty$, we change this
definition in the obvious way:
\begin{equation}\label{eq.norm.Wkinf}
    \|u\|_{\maX, \Wkinf} : = \max \|X_{1} X_{2} \ldots X_{l} u
    \|_{L^\infty}\,,
\end{equation}
the maximum being taken over the same family of vector fields.
\end{definition}

In particular,
\begin{equation}
    W^{k, p}(\Omega_0) = \{u \in L^p(\Omega_0),\ Pu \in
    L^p(\Omega_0), \text{ for all } P \in \DiffV{k}(M) \}
\end{equation}

Sometimes, when we want to stress the Lie structure $\VV$ on $M$,
we shall write $W^{k,p}(\Omega_0; M, \VV) := W^{k, p}(\Omega_0)$.

\begin{example}\label{ex.poly}\
Let $\PP$ be a curvilinear polygonal domain in the plane and let
$\Sigma(\PP)^{db}$ be the ``double'' of $\Sigma(\PP)$, which is a
Lie manifold without boundary (see Subsection \ref{ssec.poly}).
Then $\PP$ identifies with a regular open subset of
$\Sigma(\PP)^{db}$, and we have
\begin{equation*}
    \maK^m_{1}(\PP) = W^{m, 2}(\PP)
    = W^{m, 2}(\PP; \Sigma(\PP)^{db}, \VV(\PP)).
\end{equation*}
\end{example}

The following proposition shows that the second definition yields
equivalent norms.

\begin{proposition}\label{prop.Banach}\
The norms $\|\;\cdot\;\|_{\maX, \Wkp}$  and
$\|\;\cdot\;\|_{\nabla, \Wkp}$ are equivalent for any choice of
the compatible metric $g$ on $M_0$ and any choice of a system of
the finite set $\maX$ such that $\CI(M)\maX = \VV$. The spaces
$W^{k,p}(\Omega_0)$ are complete Banach spaces in the resulting
topology. Moreover, $H^{k}(\Omega_0) := W^{k,2}(\Omega_0)$ is a
Hilbert space.
\end{proposition}

\begin{proof}\
As all compatible metrics $g$ are bi-Lipschitz to each others, the
equivalence classes of the $\|\;\cdot\;\|_{\maX, \Wkp}$-norms are
independent of the choice of $g$. We will show that for any choice
$\maX$ and $g$, $\|\;\cdot\;\|_{\maX, \Wkp}$  and
$\|\;\cdot\;\|_{\nabla, \Wkp}$ are equivalent. It is clear that
then the equivalence class of $\|\;\cdot\;\|_{\maX, \Wkp}$ is
independent of the choice of $\maX$, and the equivalence class of
$\|\;\cdot\;\|_{\nabla, \Wkp}$ is independent of the choice of
$g$.

We argue by induction in $k$. The equivalence is clear for $k=0$.
We assume now that the $W^{l,p}$-norms are already equivalent for
$l=0,\dots,k-1$. Observe that if $X,Y\in \VV$, then the Koszul
formula implies $\nabla_XY\in\VV$ \cite{aln1}. To simplify
notation, we define inductively $\maX^{0}:=\maX$, and
$\maX^{i+1}=\maX^i\cup \{\nabla_X Y\,|\, X,Y\in \maX^i\}$.

By definition any $V\in \Gamma(M; T^*M^{\otimes k})$ satisfies
$(\nabla \nabla V)(X,Y)=\nabla_X \nabla_Y V - \nabla_{\nabla_X
Y}V.$ This implies for $X_1,\ldots, X_k \in \maX$
\begin{equation*}
    (\underbrace{\nabla\dots\nabla f}_{k \mbox{-times}})
    (X_1, \dots, X_k) = X_1\dots X_k f + \sum_{l=0}^{k-1}
    \sum_{Y_j\in \maX^{k-l}}
    a_{Y_1,\ldots,Y_l} Y_1\ldots Y_l\,f,
\end{equation*}
for appropriate choices of $a_{Y_1,\ldots,Y_l}\in \NN\cup\{0\}$. Hence,
\begin{equation*}
    \|(\underbrace{\nabla\dots\nabla
    f}_{k\mbox{-times}})\|_{L^p}\leq
    C \sum \|\nabla\dots\nabla f(X_1,\ldots,X_k)\|_{L^p} \leq
    C \|f\|_{\maX,W^{k,p}}.
\end{equation*}
By induction, we know that $\|Y_1,\ldots,Y_l f \|_{L^p}\leq C
\|f\|_{\nabla,W^{l,p}}$ for $Y_i\in \maX^{k-l}$, $0\leq l\leq
k-1$, and hence
\begin{multline*}
    \|X_1\ldots X_k f\|_{L^p}\leq
    \underbrace{\|\nabla\dots\nabla f\|_{L^p}
    \|X_1\|_{L^\infty} \cdots
    \|X_k\|_{L^\infty}}_{\leq C \|f\|_{\nabla, W^{k,p}}} \\
    + \underbrace{\sum_{l=0}^{k-1}\sum_{Y_1,\ldots,Y_l\in
    \maX^{k-l}} a_{Y_1,\ldots,Y_l} Y_1\ldots Y_l \, f}_{
    \leq C\|f\|_{\nabla,W^{k-1},p}}.
\end{multline*}
This implies the equivalence of the norms.

The proof of completeness is standard, see for example
\cite{Evans, Taylor3}.
\end{proof}

We shall also need the following simple observation.

\begin{lemma}\label{lemma.res}\
Let $\Omega ' \subset \Omega \subset M$ be open subsets, $\Omega_0 =
\Omega \cap M_0$, and $\Omega_0' = \Omega' \cap M_0$, $\Omega' \neq
\emptyset$. The restriction then defines continuous operators
$W^{s,p}(\Omega_0) \to W^{s, p}(\Omega'_0)$. If the various choices
($\maX,g,x_j$) are done in the same way on $\Omega$ and $\Omega'$,
then the restriction operator has norm $1$.
\end{lemma}

\subsection{Definition of Sobolev spaces using partitions of unity}
Yet another description of the spaces $W^{k,p}(\Omega_0)$ can be
obtained by using suitable partitions of unity as in \cite[Lemma
1.3]{Shubin}, whose definition we now recall. See also \cite{CGT,
Grubb, SchroheSI, SchroheFC, Skrz, Triebel}.

\begin{lemma}\label{lemma.Shubin}\ For any
$0 < \epsilon < r_{\mathrm{inj}}(M_0)/6$ there is a sequence of
points $\{x_j\} \subset M_0$, and a partition of unity $\phi_j \in
\CIc(M_0)$, such that, for some $N$ large enough depending only
on the dimension of $M_0$), we have
\begin{enumerate}[\rm (i)]
\item $\supp(\phi_j) \subset B(x_j, 2\epsilon)$;
\item $\|\nabla^k \phi_j\|_{L^\infty(M_0)}\le C_{k,\epsilon}$,
with $C_{k,\epsilon}$ independent of $j$; and
\item the sets $B(x_j, \epsilon/N)$ are disjoint, the
sets $B(x_j, \epsilon)$ form a covering of $M_0$, and the sets
$B(x_j, 4\epsilon)$ form a covering of $M_0$ of finite
multiplicity, \ie $$\sup_{y\in M_0}\# \{x_j\,|\,y\in B(x_j,
4\epsilon)\}<\infty.$$
\end{enumerate}
\end{lemma}

Fix $\epsilon\in(0,r_{\mathrm{inj}}(M_0)/6)$. Let
$\psi_j:B(x_j,4\epsilon) \to B_{\RR^n}(0,4\epsilon)$ normal
coordinates around $x_j$ (defined using the exponential map
$\exp_{x_j}:T_{x_j} M_0 \to M_0$). The uniform bounds on the
Riemann tensor $R$ and its derivatives $\nabla^k R$ imply uniform
bounds on $\nabla^k d\exp_{x_j}$, which directly implies that all
derivatives of $\psi_j$ are uniformly bounded.

\begin{proposition}\label{prop.alt.desc}\
Let $\phi_i$ and $\psi_i$ be as in the two paragraphs above. Let
$U_j=\psi_j(\Omega_0 \cap B(x_j,2\epsilon)) \subset\RR^n$. We
define
\begin{equation*}
    \nu_{k,\infty}(u) := \sup_j \|(\phi_j u)
    \circ \psi_j^{-1}\|_{W^{k,\infty}(U_j)}
\end{equation*}
and, for $1\leq p<\infty$,
\begin{equation*}
    \nu_{k,p}(u)^p := \sum_j \|(\phi_j u) \circ
    \psi_j^{-1}\|_{W^{k,p}(U_j)}^p.
\end{equation*}
Then $u \in W^{k, p}(\Omega_0)$ if, and only if, $\nu_{k,p}(u) <
\infty$. Moreover, $\nu_{k, p}(u)$ defines an equivalent norm on
$W^{k, p}(\Omega_0)$.
\end{proposition}

\begin{proof}\ We shall assume $p < \infty$, for simplicity of
notation. The case $p = \infty$ is completely similar. Consider
then $\mu(u)^p = \sum_j \|\phi_j u \|_{W^{k,p}(\Omega_0)}^p$. Then
there exists $C_{k,\varepsilon} > 0$ such that
\begin{equation}\label{eq.mu.equiv}
    C_{k,\varepsilon}^{-1}\|u\|_{W^{k, p}(\Omega_0)} \le
    \mu(u) \le C_{k,\varepsilon}\|u\|_{W^{k, p}(\Omega_0)},
\end{equation}
for all $u \in W^{k, p}(\Omega_0)$, by Lemma \ref{lemma.Shubin}
(\ie the norms are equivalent). The fact that all derivatives of
$\exp_{x_j}$ are bounded uniformly in $j$ further shows that $\mu$
and $\nu_{k,p}$ are also equivalent.
\end{proof}

The proposition gives rise to a third, equivalent definition of
Sobolev spaces. This definition is similar to the ones in
\cite{Shubin, Skrz, Triebel, Triebel2} and can be used to define the
spaces $W^{s, p}(\Omega_0)$, for any $s \in \RR$, $1 < p < \infty$,
and $\Omega_0=M_0$. The cases $p = 1$ or $p = \infty$ are more
delicate and we shall not discuss them here.

Recall that the spaces $W^{s, p}(\RR^n)$, $s \in \RR$, $ 1 < p <
\infty$ are defined using the powers of $1 + \Delta$, see
\cite[Chapter V]{Stein} or \cite[Section 13.6]{Taylor3}.

\begin{definition}[Partition of unity
definition of Sobolev spaces]\label{def.sreal}
Let $s\in \RR$, and $1<p<\infty$. Then we define
\begin{equation}\label{sreal}
    \|u\|_{W^{s,p}(M_0)}^p:=\sum_j\|(\phi_j u) \circ
\psi_j^{-1}\|_{W^{s,p}(\RR^n)}^p,\,\,\,\,\,1<p<\infty.
\end{equation}
\end{definition}

By Proposition \ref{prop.alt.desc}, this norm is equivalent to our
previous norm on $W^{s,p}(M_0)$ when $s$ is a nonnegative integer.

\begin{proposition}\label{prop.dense}\
The space $\CIc(M_0)$ is dense in $W^{s,p}(M_0)$, for $1 <
p < \infty$ and $s\in\RR$, or $1\leq p<\infty$ and
$s\in\NN\cup\{0\}$.
\end{proposition}

\begin{proof}\
For $s \in \NN\cup\{0\}$, the result is true for any manifold with
bounded geometry, see \cite[Theorem 2]{Aubin} or \cite[Theorem
2.8]{Hebey2}, or \cite{Hebey1}. For $\Omega_0 = M_0$, $s \in \RR$,
and $1 < p < \infty$, the definition of the norm on $W^{s,
p}(M_0)$ allows us to reduce right away the proof to the case of
$\RR^n$, by ignoring enough terms in the sum defining the norm
\eqref{sreal}. (We also use a cut-off function $0 \le \chi \le 1$,
$\chi \in \CIc(B_{\RR^n}(0, 4\epsilon))$, $\chi = 1$ on
$B_{\RR^n}(0, 4\epsilon)$.)
\end{proof}

We now give a characterization of the spaces $W^{s, p}(M_0)$ using
interpolation, $s \in \RR$. Let $k \in \NN \cup
\{0\}$ and let $\widetilde W^{-k, p}(M_0)$ be the set of
distributions on $M_0$ that extend by continuity to linear
functionals on $W^{k, q}(M_0)$, $p^{-1} + q^{-1} = 1$, using
Proposition \ref{prop.dense}. That is, let $\widetilde W^{-k,
p}(M_0)$ be the set of distributions on $M_0$ that define
continuous linear functionals on $W^{k, q}(M_0)$, $p^{-1} + q^{-1}
= 1$. We let
\begin{equation*}
    \widetilde W^{\theta k, k, p}(M_0) :=
    [\widetilde W^{0, p}(M_0), W^{k, p}(M_0)]_{\theta}\,,
    \quad 0 \le \theta \le 1\,,
\end{equation*}
be the complex interpolation spaces. Similarly, we define
\begin{equation*}
    \widetilde W^{-\theta k, k, p}(M_0) = [\widetilde W^{0,
    p}(M_0), W^{-k,p}(M_0)]_{\theta}.
\end{equation*}
(See \cite{BerghLofstrom} or \cite[Chapter 4]{Taylor1} for the
definition of the complex interpolation spaces.)

The following proposition is an analogue of Proposition
\ref{prop.alt.desc}. Its main role is to give an intrinsic
definition of the spaces $W^{s, p}(M_0)$, a definition that is
independent of choices.

\begin{proposition}\label{prop.alt.desc1}\
Let $1 < p < \infty$ and $k > |s|$. Then we have a topological
equality $\widetilde W^{s, k , p}(M_0) = W^{s, p}(M_0)$. In
particular, the spaces $W^{s, p}(M_0)$, $s \in \RR$, do not depend
on the choice of the covering $B(x_j, \epsilon)$ and of the
subordinated partition of unity and we have
\begin{equation*}
    [W^{s, p}(M_0), W^{0 , p}(M_0)]_{\theta} = W^{\theta s,
    p}(M_0)\,,\quad 0 \le \theta \le 1\,.
\end{equation*}
Moreover, the pairing between functions and distributions defines
an isomorphism $W^{s, p}(M_0)^* \simeq W^{-s, q}(M_0)$, where $1/p
+ 1/q = 1$.
\end{proposition}

\begin{proof}\
This proposition is known if $M_0 = \RR^n$ with the usual metric
\cite{Taylor3}[Equation (6.5), page 23]. In particular, 
$\widetilde W^{s,p}(\RR^n) = W^{s, p}(\RR^n)$. As in the proof of Proposition
\ref{prop.alt.desc} one shows that the quantity
\begin{equation}\label{eq.nu.sp}
    \nu_{s, p}(u)^p := \sum_j\|(\phi_j u) \circ
    \psi_j^{-1}\|_{\tilde W^{s,p}(\RR^n)}^p,
\end{equation}
is equivalent to the norm on $\tilde W^{s, p}(M_0)$. This implies
$\tilde W^{s,p}(M_0) = W^{s, p}(M_0)$.

Choose $k$ large. Then we have
\begin{multline*}
    [W^{s, p}(M_0), W^{0, p}(M_0)]_{\theta} = [W^{s, k, p}(M_0), W^{0,
    k, p}(M_0)]_{\theta}\\ = W^{\theta s, k, p}(M_0) = W^{\theta s,
    p}(M_0).
\end{multline*}
The last part follows from the compatibility of interpolation with
taking duals. This completes the proof.
\end{proof}

The above proposition provides us with several corollaries. First,
from the interpolation properties of the spaces $W^{s, p}(M_0)$, we
obtain the following corollary.

\begin{corollary}\label{cor.n.mult}\
Let $\phi \in W^{k, \infty}(M_0)$, $k\in \NN\cup\{0\}$, $p \in (1,
\infty)$, and $s\in\RR$  with $k \ge |s|$. Then multiplication by
$\phi$ defines a bounded operator on $W^{s,p}(M_0)$ of norm at
most $C_k\|\phi\|_{W^{k, \infty}(M_0)}$. Similarly, any
differential operator $P \in \DiffV{m}(M)$ defines continuous maps
$P : W^{s, p}(M_0) \to W^{s-m, p}(M_0)$.
\end{corollary}

\begin{proof}\
For $s \in \NN\cup\{0\}$, this follows from the definition of the norm on
$W^{k, \infty}(M_0)$ and from the definition of $\DiffV{m}(M)$ as the
linear span of differential operators of the form $fX_1 \ldots X_k$,
($f \in \CI(M) \subset W^{k, \infty}$, $X_j \in \VV$, and $0
\le k \le m$), and from the definition of the spaces
$W^{k,p}(\Omega_0)$.

For $s \le m$, the statement follows by duality. For the other values of $s$,
the result follows by interpolation.
\end{proof}

Next, recall that an isomorphism $\phi: M \to M'$ of the Lie
manifolds $(M_0, M, \VV)$ and $(M_0', M', \VV')$ is defined to be a
diffeomorphism such that $\phi_*(\VV) = \VV'$. We then have the
following invariance property of the Sobolev spaces that we have
introduced.

\begin{corollary}\label{cor.invariance}\
Let $\phi : M \to M'$ be an isomorphism of Lie manifolds,
$\Omega_0 \subset M_0$ be an open subset and $\Omega' =
\phi(\Omega)$. Let $p \in [1, \infty]$ if $s \in \NN\cup\{0\}$, and $p
\in (1, \infty)$ if $s \not \in \NN\cup\{0\}$. Then $f \to f \circ \phi$
extends to an isomorphism $\phi^* : W^{s, p}(\Omega') \to W^{s,
p}(\Omega)$ of Banach spaces.
\end{corollary}

\begin{proof}\
For $s \in \NN\cup\{0\}$, this follows right away from definitions and
Proposition~\ref{prop.Banach}. For $-s \in \NN\cup\{0\}$, this
follows by duality, Proposition \eqref{prop.alt.desc1}. For the other
values of $s$, the result follows from the same proposition, by
interpolation.
\end{proof}

Recall now that $M_0$ is complete \cite{aln1}. Hence the Laplace
operator $\Delta = \nabla^*\nabla$ is essentially self-adjoint on
$\CIc(M_{0})$ by \cite{gaffney, roel}. We shall define then $(1 +
\Delta)^{s/2}$ using the spectral theorem.

\begin{proposition}\label{prop.domain}\ The space
$H^s(M_0) := W^{s, 2}(M_0)$, $s \ge 0$, identifies with the domain
of $(1 + \Delta)^{s/2}$, if we endow the latter with the graph
topology.
\end{proposition}

\begin{proof}\ For $s \in \NN\cup\{0\}$, the result is true for
any manifold of bounded geometry, by \cite[Proposition 3]{Aubin}.
For $s \in \RR$, the result follows from interpolation, because
the interpolation spaces are compatible with powers of operators
(see, for example, the chapter on Sobolev spaces in Taylor's book
\cite{Taylor1}).
\end{proof}

The well known Gagliardo--Nirenberg--Sobolev inequality
\cite{Aubin, Evans, Hebey2} holds also in our
setting.

\begin{proposition}\label{prop.Sobolev}\ Denote by $n$ the dimension
of $M_0$. Assume that $1/p = 1/q - m/n$, $1 < q \leq p < \infty$,
where $m\geq 0$. Then $W^{s, q}(M_0)$ is continuously embedded in
$W^{s-m, p}(M_0)$.
\end{proposition}

\begin{proof}\
If $s$ and $m$ are integers, $s\geq m\geq 0$, the statement of the
proposition is true for manifolds with bounded geometry,
\cite[Theorem 7]{Aubin} or \cite[Corollary 3.1.9]{Hebey2}. By
duality (see Proposition \ref{prop.alt.desc1}), we obtain the same
result when $s \le 0$, $s\in \ZZ$. Then, for integer $s,m$,
$0<s<m$ we obtain the corresponding embedding by composition
$W^{s,q}(M_0)\to W^{0,r}(M_0)\to W^{s-m,p}(M_0)$, with
$1/r=1/q-s/n$. This proves the result for integral values of $s$.
For non-integral values of $s$, the result follows by
interpolation using again Proposition~\ref{prop.alt.desc1}.
\end{proof}

The Rellich-Kondrachov's theorem on the compactness of the
embeddings of Proposition \ref{prop.Sobolev} for $1/p > 1/q -m/n$
is true if $M_0$ is compact \cite[Theorem 9]{Aubin}. This happens
precisely when $M = M_0$, which is a trivial case of a manifold
with a Lie structure at infinity. On the other hand, it is easily
seen (and well known) that this compactness cannot be true for
$M_0$ non-compact. We will nevertheless obtain compactness in the
next section by using Sobolev spaces with weights, see
Theorem~\ref{theorem.Kondrachov}.

\section{Sobolev spaces on regular open subsets\label{sec.MB}}

Let $\Omega \subset M$ be an open subset. Recall that $\Omega$ is
a regular open subset in $M$ if, and only if, $\Omega$ and
$\overline{\Omega}$ have the same boundary in $M$, denoted
$\patop\overline\Omega$,, and if $\patop\overline\Omega$ is a
regular submanifold of $M$. Let $\Omega_0 = \Omega \cap M_0$. Then
$\pa \Omega_0 := (\pa \Omega) \cap M_0= \patop \overline\Omega
\cap M_0$ is a smooth submanifold of codimension one of $M_0$ (see
Figure \ref{fig.sigma}). We shall denote $W^{s,
p}(\overline{\Omega}) = W^{s, p}(\Omega) = W^{s, p}(\Omega_0)$.
Throughout this section $\Omega$ will denote a regular open
subset of $M$.

We have the following analogue of the classical extension theorem.

\begin{theorem}\label{theorem.ext}\
Let $\Omega \subset M$ be a regular open subset. Then
there exists a linear operator $E$ mapping measurable functions on
$\Omega_0$ to measurable functions on $M_0$ with the properties:
\begin{enumerate}[\rm (i)]
\item\ $E$ maps $W^{k,p}(\Omega_0)$ continuously into $\WkpM$ for
every $p\in[1,\infty]$ and every integer $k\geq 0$, and
\item\ $E u \vert_{\Omega_0} = u$.
\end{enumerate}
\end{theorem}

\begin{proof}\
Since $\pa\Omega_0$ is a regular submanifold we can fix a compatible
metric $g$ on $M_0$ and a tubular neighborhood $V_0$ of $\pa\Omega_0$
such that $V_0 \simeq (\pa\Omega_0) \times (-\varepsilon_0,
\varepsilon_0)$, $\varepsilon_0>0$. Let
$\varepsilon=\min(\varepsilon_0,r_{\mathrm{inj}}(M_0))/20$, where
$r_{\mathrm{inj}}(M_0)>0$ is the injectivity radius of $M_0$. By
Zorn's lemma and the fact that $M_0$ has bounded geometry we can
choose a maximal, countable set of disjoint balls
$B(x_i,\varepsilon)$, $i\in I$. Since this family of balls is maximal
we have $M_0=\cup_iB(x_i,2\varepsilon)$. For each $i$ we fix a smooth
function $\eta_i$ supported in $B(x_i,3\varepsilon)$ and equal to $1$
in $B(x_i,2\varepsilon)$. This can be done easily in local coordinates
around the point $x_i$; since the metric $g$ is induced by a metric
$g$ on $A$ we may also assume that all derivatives of order up to $k$
of $\eta_i$ are bounded by a constant $C_{k,\varepsilon}$ independent
of $i$.  We then set $\tilde\eta_i := \big(\sum_{j\in I} \eta_j^2
\big)^{-1/2} \eta_i.$ Then $\sum_{i\in I} \tilde\eta_i^2 = 1$,
$\tilde\eta_i$ equals $1$ on $B(x_i,\ep)$ and is supported in
$B(x_i,3\ep)$.

Following \cite[Ch. 6]{Stein} we also define two smooth cutoff
functions adapted to the set $\Omega_0$. We start with a function
$\psi:\RR\to[0,1]$ which is equal to $1$ on $[-3,3]$ and which has
support in $[-6,6]$

Let $\varphi = (\varphi_1,\varphi_2)$ denote the isomorphism
between $V_0$ and $\pa\Omega_0 \times
(-\varepsilon_0,\varepsilon_0)$, where $\varphi_1 : V_0 \to
\pa\Omega_0$ and $\varphi_2 :
V_0\to(-\varepsilon_0,\varepsilon_0)$. We define
\begin{equation*}
\Lambda_+(x):=
\begin{cases}
    \ \ 0 & \text{if }\; x\in M_0\setminus V_0\\ \
    \psi(\varphi_2(x)/\varepsilon) & \text{if }\; x\in V_0,
\end{cases}
\end{equation*}
and $\Lambda_-(x):=1-\Lambda_+(x)$. Clearly $\Lambda_+$ and
$\Lambda_-$ are smooth functions on $M_0$ and $\Lambda_+(x) +
\Lambda_-(x)= 1$. Obviously, $\Lambda_+$ is supported in a
neighborhood of $\pa\Omega_0$ and $\Lambda_-$ is supported in the
complement of a neighborhood of $\pa\Omega_0$.

Let $\pa\Omega_0 = A_1 \cup A_2 \cup \ldots$ denote the
decomposition of $\pa\Omega_0$ into connected components. Let $V_0
= B_1 \cup B_2 \cup \ldots$ denote the corresponding decomposition
of $V_0$ into connected components, namely, $B_j =
\varphi^{-1}(A_j \times (-\varepsilon_0, \varepsilon_0))$. Since
$\pa\Omega_0=\pa{\overline{\Omega}_0}$, we have
$\varphi(\Omega_0\cap B_j)=A_j\times(-\varepsilon_0,0)$ or
$\varphi(\Omega_0\cap B_j)=A_j\times(0,\varepsilon_0)$. Thus, if
necessary, we may change the sign of $\varphi$ on some of the
connected components of $V_0$ in such a way that
\begin{equation*}
    \varphi(\Omega_0\cap V_0) = \pa\Omega_0 \times
    (0,\varepsilon_0).
\end{equation*}
Let $\psi_0$ denote a fixed smooth function, $\psi_0 :
\RR\to[0,1]$, $\psi_0(t)=1$ if $t\geq-\varepsilon$ and
$\psi_0(t)=0$ if $t\leq-2\varepsilon$, and let
\begin{equation*}
\Lambda_0(x)=
\begin{cases}
\ \ 1 & \text{if }\; x \in \Omega_0 \setminus V_0\\ \ \ 0 &
\text{if }\; x \in M_0 \setminus (\Omega_0 \cup V_0) \\ \
\psi_0(\varphi_2(x)) & \text{if }\; x \in V_0.
\end{cases}
\end{equation*}

We look now at the points $x_i$ defined in the first paragraph of
the proof. Let $J_1=\{i\in I:d(x_i,\pa\Omega_0)\leq 10
\varepsilon\}$ and $J_2=\{i\in I : d(x_i,\pa\Omega_0)
> 10\varepsilon\}$. For every point $x_i$, $i\in J_1$, there is a
point $y_i\in\pa\Omega_0$ with the property that
$B(x_i,4\varepsilon) \subset B(y_i,15\varepsilon)$. Let
$B_{\pa\Omega_0}(y_i, 15\varepsilon)$ denote the ball in
$\pa\Omega_0$ of center $y_i$ and radius $15\varepsilon$ (with
respect to the induced metric on $\pa\Omega_0$). Let $h_i :
B_{\pa\Omega_0}(y_i,15\varepsilon) \to B_{\RR^{n-1}}(0,
15\varepsilon)$ denote the normal system of coordinates around the
point $y_i$. Finally let $g_i : B_{\RR^{n-1}}(0, 15\varepsilon)
\times (-15\varepsilon, 15\varepsilon) \to V_0$ denote the map
$g_i(v,t) = \varphi^{-1}(h_i^{-1}(v),t)$.

Let $E_{\RR^n}$ denote the extension operator that maps
$W^{k,p}(\RR^n_+)$ to $W^{k,p}(\RR^n)$ continuously, where
$\RR^n_+$ denotes the half-space $\{x:x_n > 0\}$. Clearly,
$E_{\RR^n}u\vert_{\RR^n_+} = u$. The existence of this extension
operator is a classical fact, for instance, see \cite[Chapter~6]{Stein}.
For any $u\in W^{k,p}(\Omega_0)$ and $i\in J_1$ the
function $(\tilde\eta_iu)\circ g_i$ is well defined on $\RR^n_+$ simply
by setting it equal to $0$ outside the set
$B_{\RR^{n-1}}(0,15\varepsilon)\times(0, 15\varepsilon)$. Clearly,
$(\tilde\eta_i u)\circ g_i\in W^{k,p}(\RR^n_+)$. We define the extension
$Eu$ by the formula
\begin{equation}\label{extension}
    Eu(x) = \Lambda_0(x) \Lambda_-(x)u(x) + \Lambda_0(x)
    \Lambda_+(x) \sum_{i\in J_1} \tilde\eta_i(x)
   \Bigl(E_{\RR^n}[(\tilde\eta_iu)\circ g_i]\Bigr)(g_i^{-1}x)\,.
\end{equation}
Notice that for all $i\in J_2$, the function $\tilde\eta_i$
vanishes on the support of $\Lambda_+$, and hence
\begin{equation}\label{verify}
    \sum_{i\in J_1}\tilde\eta_i^2(x)
    = \sum_{i\in I}\tilde\eta_i^2(x) = 1 \;
    \text{ in }\mathrm{supp}\, \Lambda_+.
\end{equation}
This formula implies $Eu\vert_{\Omega_0}=u$. It remains to verify
that
\begin{equation*}
    \|Eu\|_{W^{k,p}(M_0)}\leq C_k\|u\|_{W^{k,p}(\Omega_0)}.
\end{equation*}
This follows as in \cite{Stein} using \eqref{verify}, the fact
that the extension $E_{\RR^n}$ satisfies the same bound, and the
definition of the Sobolev spaces using partitions of unity
(Proposition \ref{prop.alt.desc}).
\end{proof}

Let $\Omega$ be a regular open subset of $M$ and $\Omega_0 =
\Omega \cap M$, as before. We shall denote by
$\overline{\Omega}_0$ the closure of $\Omega_0$ in $M_0$.

\begin{theorem}\label{theorem.dense}\
The space $\CIc(\overline{\Omega}_0)$ is dense in
$W^{k,p}(\Omega_0)$, for $1 \le p < \infty$.
\end{theorem}

\begin{proof}\ For any $u\in W^{k,p}(\Omega_0)$ let $Eu$
denote its extension from Theorem \ref{theorem.ext}, $Eu\in
W^{k,p}(M_0)$. By Proposition \ref{prop.dense}, there is a
sequence of functions $f_j\in C_c^{\infty}(M_0)$ with the property
that
\begin{equation*}
\lim_{j\to\infty}f_j=Eu\text{ in }W^{k,p}(M_0).
\end{equation*}
Thus $\lim_{j\to\infty}f_j\vert_{{\Omega_0}}=u$ in
$W^{k,p}(\Omega_0)$, as desired.
\end{proof}

\begin{theorem}\label{theorem.trace.I}\
The restriction map $\CIc(\overline{\Omega}_0) \to
\CIc(\pa\Omega_0)$ extends to a continuous map $T :
W^{k,p}(\Omega_0) \to W^{k-1, p}(\pa\Omega_0)$, for $1 \le p \leq
\infty$.
\end{theorem}

\begin{proof}\ 
The case $p=\infty$ is obvious.  In the case $1 \leq p < \infty$, we
shall assume that the compatible metric on $M_0$ restricts to a
product type metric on $V_0$, our distinguished tubular neighborhood
of $\pa\Omega_0$. As the curvature of $M_0$ and the second fundamental
form of $\pa\Omega_0$ in $M_0$ are bounded (see
Corollary~\ref{cor.II.bdd}), there is an $\ep_1>0$ such that, in
normal coordinates, the hypersurface $\pa\Omega_0$ is the graph of a
function on balls of radius $\leq \ep_1$.

We use the definitions of the Sobolev spaces using partitions of
unity, Proposition \ref{prop.alt.desc} and Lemma
\ref{lemma.Shubin} with $\varepsilon = \min(\epsilon_1,\epsilon_0,
r_{\mathrm{inj}}(M_0))/10$. Let $B(x_j, 2\varepsilon)$ denote the
balls in the cover of $M_0$ in Lemma \ref{lemma.Shubin}, let
$\psi_j:B(\ep,x_j)\to B(\ep,0)$ denote normal coordinates based in
$x_j$, and let $1=\sum_j\phi_j$ be a corresponding partition of
unity. Then $\widetilde{\phi}_j = \phi_j\vert_{\pa\Omega_0}$ form
a partition of unity on $\pa\Omega_0$.

Start with a function $u\in W^{k,p}(\Omega_0)$ and let
$u_j=(u\phi_j)\circ\psi_j^{-1}$, $u_j\in
W^{k,p}(\psi_j(\Omega_0\cap B(x_j,4\varepsilon)))$. In addition
$u_j\equiv 0$ outside the set $\psi_j(\Omega_0\cap
B(x_j,2\varepsilon))$. If
$B(x_j,4\varepsilon)\cap\pa\Omega_0=\emptyset$ let
$\widetilde{T}(u_j)=0$. Otherwise notice that
$B(x_j,4\varepsilon)$ is included in $V_0$, the tubular
neighborhood of $\pa\Omega_0$, thus the set
$\psi_j(\pa\Omega_0\cap B(x_j,4\varepsilon))$ is the intersection
of a graph and the ball $B_{\RR^n}(0,4\varepsilon)$. We can then
let $\widetilde{T}(u_j)$ denote the Euclidean restriction of $u_j$
to $\psi_j(\pa\Omega_0\cap B(x_j,4\varepsilon))$ (see
\cite[Section 5.5]{Evans}). Clearly $\widetilde{T}(u_j)$ is
supported in $\psi_j(\pa\Omega_0\cap B(x_j,2\varepsilon))$ and
\begin{equation*}
    \|\widetilde{T}(u_j) \circ
    \widetilde{\psi}_j\|_{W^{k-1,p}(\pa\Omega_0)}
    \leq C\|u_j\|_{W^{k,p}(\psi_j(\Omega_0\cap
    B(x_j,4\varepsilon)))},
\end{equation*}
where $\widetilde{\psi}_j=\psi_j\vert_{\Omega_0}$ and the constant
$C$ is independent of $j$ (recall that $\psi_j(\pa\Omega_0\cap
B(x_j,4\varepsilon))$ is the intersection of a hyperplane and the
ball $B_{\RR^n}(0,4\varepsilon)$). Let
\begin{equation*}
    Tu=\sum_j\widetilde{T}(u_j)\circ\widetilde{\psi}_j.
\end{equation*}
Since the sum is uniformly locally finite, $Tu$ is well-defined
and we have
\begin{equation*}
\begin{split}
    \|Tu\|_{W^{k-1,p}(\pa\Omega_0)}^p&\leq
    C\sum_j\|\widetilde{T}(u_j)\circ
    \widetilde{\psi}_j\|_{W^{k-1,p}(\pa\Omega_0)}^p\\
    &\leq C\sum_j\|u_j\|_{W^{k,p}(\psi_j(\Omega_0\cap
    B(x_j,4\varepsilon)))}^p\leq C\|u\|_{W^{k,p(\Omega_0)}}^p,
\end{split}
\end{equation*}
with constants $C$ independent of $u$. The fact that
$Tu\vert_{\mathcal{C}_c^\infty(\Omega_0)}$ is indeed the
restriction operator follows immediately from the definition.
\end{proof}

We shall see that if $p=2$, we get a surjective map $W^{s,
2}(\Omega_0) \to W^{s-1/2, 2}(\pa\Omega_0)$ (Theorem
\ref{theorem.res}).

In the following, $\pa_\nu$ denotes derivative in the normal
direction of the hypersurface $\pa \Omega_0\subset M_0$.

\begin{theorem}\label{theorem.zero}\
The closure of $\CIc(\Omega_0)$ in $\WkpO$ is the intersection of the
kernels of $T \circ \pa_\nu^j:\WkpO\to W^{k-j-1,p}(\Omega_0)$, $0 \le
j \le k-1$, $1 \le p < \infty$.
\end{theorem}

\begin{proof}\
The proof is reduced to the Euclidean case \cite{Adams,
Evans,LionsMagenes1, Taylor1} following the same pattern of reasoning
as in the previous theorem.
\end{proof}

The Gagliardo--Nirenberg--Sobolev theorem holds also for manifolds
with boundary.

\begin{theorem}\label{theorem.Sobolev}\
Denote by $n$ the dimension of $M$ and let $\Omega \subset M$ be a
regular open subset in $M$. Assume that $1/p = 1/q - m/n > 0$, $1
\le q < \infty$, where $m \le k$ is an integer. Then $W^{k,
q}(\Omega_0)$ is continuously embedded in $W^{k-m, p}(\Omega_0)$.
\end{theorem}

\begin{proof}\  This can be proved using Proposition
\ref{prop.Sobolev} and Theorem \ref{theorem.ext}. Indeed, denote
by
\begin{equation*}
    j : W^{k, q}(M_0)  \to  W^{k-m, p}(M_0)
\end{equation*}
the continuous inclusion of Proposition \ref{prop.Sobolev}. Also,
denote by $r$ the restriction maps $\WkpM \to \WkpO$. Then the
maps
\begin{equation*}
    W^{k, q}(\Omega_0) \stackrel{E}{\longrightarrow} W^{k, q}(M_0)
    \stackrel{j}{\longrightarrow} W^{k-m, p}(M_0)
    \stackrel{r}{\longrightarrow} W^{k-m, p}(\Omega_0)
\end{equation*}
are well defined and continuous. Their composition is the
inclusion of $W^{k, q}(\Omega_0)$ into $W^{k-m, p}(\Omega_0)$.
This completes the proof.
\end{proof}

For the proof of a variant of Rellich--Kondrachov's compactness
theorem, we shall need Sobolev spaces with weights. Let $\Omega
\subset M$ be a regular open subset. Let $a_H \in \RR$ be a
parameter associated to each boundary hyperface (\ie face of
codimension one) of the manifold with corners $\overline{\Omega}$.
Fix for any boundary hyperface $H \subset \overline{\Omega}$ a
defining function $\rho_H$, that is a function $\rho_H \ge 0$ such
that $H =\{\rho_H = 0\}$ and $d\rho_H \not = 0$ on $H$. Let
\begin{equation}\label{eq.def.rho}
    \rho = \prod \rho_H^{a_H},
\end{equation}
the product being taken over all boundary hyperfaces of
$\overline{\Omega}$. A~function of the form $\psi\rho$, with $\psi
> 0$, $\psi$ smooth on $\overline{\Omega}$, and $\rho$ as in
Equation \eqref{eq.def.rho} will be called an {\em admissible
weight} of $\overline{\Omega}$ (or simply an admissible weight
when $\Omega$ is understood). We define then the weighted Sobolev
space $\WkpO$ by
\begin{equation}\label{eq.def.weights}
    \rho \WkpO := \{\rho u, u \in \WkpO\},
\end{equation}
with the norm $\|\rho^su\|_{\rho^s \WkpO} := \|u\|_{\WkpO}$.

Note that in the definition of an admissible weight of
$\overline{\Omega}$, for a regular open subset $\Omega \subset M$
of the Lie manifold $(M, \VV)$, we allow also powers of the
defining functions of the boundary hyperfaces contained in $\pa
\Omega = \patop \overline{\Omega}$, the true boundary of
$\overline{\Omega}$. In the next compactness theorem, however, we
shall allow only the powers of the defining functions of $M$, or,
which is the same thing, only powers of the defining functions of
the boundary hyperfaces of $\overline{\Omega}$ whose union is
$\pa_\infty \overline{\Omega}$ (see Figure \ref{fig.sigma}).

\begin{theorem}\label{theorem.Kondrachov}\
Denote by $n$ the dimension of $M$ and let $\Omega \subset M$ be a
regular open subset, $\Omega_0 = \Omega \cap M_0$. Assume that
$1/p > 1/q - m/n > 0$, $1 \le q < \infty$, where $m \in
\{1,\ldots, k\}$ is an integer, and that $s > s'$ are real
parameters. Then $\rho^sW^{k, q}(\Omega_0)$ is compactly embedded
in $\rho^{s'}W^{k-m, p}(\Omega_0)$ for any admissible weight
$\rho := \prod_H \rho_H^{a_H}$ of $M$ such that $a_H > 0$
for any boundary hyperface $H$ of $M$.
\end{theorem}

\begin{proof}\ The same argument as that in the proof of Theorem
\ref{theorem.Sobolev} allows us to assume that $\Omega_0 = M_0$.
The norms are chosen such that $\WkpO \ni u \mapsto \rho^s u \in
\rho^s\WkpO$ is an isometry. Thus, it is enough to prove that
$\rho^s : W^{k, q}(\Omega_0) \to W^{k - m, p}(\Omega_0)$, $s > 0$,
is a compact operator.

For any defining function $\rho_H$ and any $X \in \VV$, we have
that $X(\rho_H)$ vanishes on $H$, since $X$ is tangent to $H$. We
obtain that $X(\rho^s) = \rho^s f_X$, for some $f_X \in \CI(M)$.
Then, by induction, $X_1 X_2 \ldots X_k (\rho^s) = \rho^s g$, for
some $g \in \CI(M)$.

Let $\chi \in \CI([0, \infty)$ be equal to $0$ on $[0, 1/2]$,
equal to $1$ on $[1, \infty)$, and non-negative everywhere. Define
$\phi_\epsilon = \chi(\epsilon^{-1} \rho^s)$. Then
\begin{equation*}
    \| X_1 X_2 \ldots X_k \big(\rho^s \phi_\epsilon
    - \rho^s\big) \|_{L^\infty} \to 0\,, \quad \text{as }\;
    \epsilon \to 0,
\end{equation*}
for any $X_1, X_2, \ldots, X_k \in \VV$. Corollary
\ref{cor.n.mult} then shows that $\rho^s \phi_\epsilon \mapsto \rho^s$ in
the norm of bounded operators on $W^{s, p}(\Omega_0)$. But
multiplication by $\rho^s\phi_\epsilon$ is a compact operator, by the
Rellich-Kondrachov's theorem for compact manifolds with boundary
\cite[Theorem 9]{Aubin}. This completes the proof.
\end{proof}

We end with the following generalization of the classical
restriction theorem for the Hilbertian Sobolev spaces $H^s(M_0) :=
W^{s, 2}(M_0)$.

\begin{theorem}\label{theorem.res}\
Let $N_0 \subset M_0$ be a tame submanifold of codimension $k$ of
the Lie manifold $(M_0, M, \VV)$. Restriction of smooth functions
extends to a bounded, surjective map
\begin{equation*}
    H^{s}(M_0) \to H^{s - k/2}(N_0),
\end{equation*}
for any $s > k/2$. In particular, $H^{s}(\Omega_0) \to H^{s -
1/2}(\pa\Omega_0)$ is continuous and surjective.
\end{theorem}

\begin{proof}\
Let $B \to N$ be the vector bundle defining the Lie structure at
infinity $(N, B)$ on $N_0$ and $A \to M$ be the vector bundle defining
the Lie structure at infinity $(M, A)$ on $M_0$. (See Section
\ref{sec.SM} for further explanation of this notation.)  The existence
of tubular neighborhoods, Theorem \ref{thm.tub.neigh}, and a partition
of unity argument, allows us to assume that $M = N \times S^1$ and
that $A = B \times TS^1$ (external product). Since the Sobolev spaces
$H^s(M_0)$ and $H^{s-1/2}(N_0)$ do not depend on the metric on $A$ and
$B$, we can assume that the circle $S^1$ is given the invariant metric
making it of length $2\pi$ and that $M_0$ is given the product metric.
The rest of the proof now is independent of the way we obtain the
product metric on $M_0$.

Let $S^1$ be the unit circle in the plane. Let us denote by
$\Delta_M, \Delta_N,$ and $\Delta_{S^1}$ the Laplace operators on
$M_0$, $N_0$, and $S^1$, respectively. Then $\Delta_{M} =
\Delta_{N} + \Delta_{S^1}$ and $\Delta_{S^1}= -\pa^2/\pa \theta^2$
has spectrum $\{4\pi^2n^2\,|\,n \in \NN\cup\{0\}\}$. We can
decompose $L^2(N_0 \times S^1)$ according to the eigenvalues $n
\in \ZZ$ of $-\frac{1}{2\pi\imath}\pa_\theta$:
\begin{equation*}
    L^2(N_0 \times S^1) \simeq \bigoplus_{n \in \ZZ}
    L^2(N_0 \times S^1)_n \simeq \bigoplus_{n \in \ZZ} L^2(N_0),
\end{equation*}
where the isomorphism $L^2(N_0 \times S^1)_n \simeq L^2(N_0)$ is
obtained by restricting to $N_0 = N_0 \times \{1\}$, $1 \in S^1$.
We use this isomorphism to identify the above spaces in what follows.

Let $\xi \in L^2(N_0 \times S^1)$. Then $\xi$ identifies
with a sequence $(\xi_n)$ under the above isomorphism. By
Proposition \ref{prop.domain}, we have that $\xi \in H^s(N_0
\times S^1)$ if, and only if, $(1 + \Delta_{M})^{s/2} \xi =
\sum_n\big((1 + n^2 + \Delta_{N})^{s/2}\xi_n \big) \in \bigoplus_{n \in
\ZZ} L^2(N_0) \simeq L^2(N_0 \times S^1)$. The restriction of
$\xi$ to $N_0$ is then given by $\sum_n \xi_n$. We want to show
that $\sum_n \xi_n \in H^{s-1/2}(N_0)$, which is equivalent to $(1
+ \Delta_{N})^{s/2-1/4}\big(\sum \xi_n\big) \in L^2(N_0)$.

The spectral spaces of $\Delta_{N}$ corresponding to $[m,m+1) \subset \RR$, 
$m\in \NN\cup\{0\}$ give an orthogonal direct sum decomposition of $L^2(N_0)$.

We decompose $\xi_n = \sum_m \xi_{mn}$,
with $\xi_{mn}$ in the spectral space corresponding to $[m, m+1)$ of
$\Delta_{N}$.  
Note that $\xi_{mn}$ is orthogonal to $\xi_{m'n}$ for $m\neq m'$.
Denote $h =(1+ m^2)^{-1/2}$, $f(t)= (1 + t^2)^{-s}$, and $C = 1 + \int_\RR
f(t)dt$. Then an application of the Cauchy--Schwartz inequality gives
\begin{multline}\label{eq.C}
    	(1 + m^2)^{s - 1/2} \Bigl( \sum_n \|\xi_{mn} \| \Bigr )^2 \\ \le
	(1 + m^2)^{s - 1/2} \Bigl(\sum_n (1 + n^2 + m^2)^{-s} \Bigr
	) \sum_n \| (1 + n^2 + m^2)^{s/2} \xi_{mn}\|^2 \\ 
	\le h \Bigl(\sum_{n}f(nh) \Bigr) \sum_n
	\| (1 + n^2 + m^2)^{s/2} \xi_{mn}\|^2
	\le C_s \sum_n
	\| (1 + n^2 + m^2)^{s/2} \xi_{mn}\|^2.
\end{multline}
The constant $C_s$ is independent of $m$ (but depends on $s$). 
We sum over $m$ and obtain
\begin{multline}\label{eq.m}
    \| \sum_n (1 + \Delta_{N})^{s/2-1/4} \xi_n \|^2
    = \sum_m \|\sum_n (1 + \Delta_{N})^{s/2-1/4} \xi_{nm}\|^2
    \\ \le \sum_m (1 + (m+1)^2)^{s - 1/2}
    \Bigl( \sum_n \| \xi_{nm} \| \Bigr )^2
    \le 2^s \sum_m (1 + m^2)^{s - 1/2}
    \Bigl( \sum_n \| \xi_{nm} \| \Bigr )^2 \\
    \le 2^sC_s \sum_{n, m} \| (1 + n^2 + m^2)^{s/2} \xi_{nm}\|^2
    \le 2^sC_s \sum_{n, m} \| (1 + n^2 + \Delta_{N})^{s/2} \xi_{nm}\|^2\\ 
    = 2^s C_s \sum_{n} \| (1 + n^2 + \Delta_{N})^{s/2} \xi_{n}\|^2 ,
\end{multline}
with the same constant $C_s$ as in Equation \eqref{eq.C}.
This shows that $\zeta := \sum_n \xi_n \in H^{s-1/2}(N_0)$ if $\xi
= (\xi_n) \in \bigoplus_n L^2(N_0) \simeq L^2(N_0 \times S^1)$ is a
{\em finite} sequence such that
    $\|\xi\|_{H^s} := \sum_n \|( 1 + n^2 + \Delta_{N})^{s/2}
    \xi_n\|^2_{L^2(N_0)} < \infty,$
and that $\zeta$ depends continuously on $\xi \in H^s(N_0 \times
S^1)$. This completes the proof.
\end{proof}

We finally obtain the following consequences for a curvilinear
polygonal domain~$\PP$ (see Subsection \ref{ssec.poly}). First,
recall that the distance $\vartheta(x)$ from $x$ to the vertices
of a curvilinear polygon $\PP$ and $r_{\PP}$ have bounded
quotients, and hence define the same weighted Sobolev spaces
(Equation~\eqref{eq.same.Sobo}). Moreover, the function $r_{\PP}$
is an admissible weight. Recall that $\PP$ has a compactification
$\Sigma(\PP)$ that is a Lie manifold with boundary (that is, the
closure of a regular open subset of a Lie manifold $M$). Let us
write $W^{m,p}(\Sigma(\PP)) := W^{m, p}(\PP)$ the Sobolev spaces
defined by the structural Lie algebra of vector fields on
$\Sigma(\PP)$. Then
\begin{equation}\label{eq.def.wssp}
    \maK_{a}^m(\PP; \vartheta) = r_{\Omega}^{a-1}
    \maK_{1}^m(\PP; r_{\PP}) =
    r_{\PP}^{a-1}W^{m, 2}(\Sigma(\PP)).
\end{equation}
This identifies the weighted Sobolev spaces on $\PP$
with a weighted Sobolev space of the form $\rho \WkpO$.

Motivated by Equation \eqref{eq.def.wssp}, we now define
\begin{equation}\label{eq.def.wsspB}
    \maK_{a}^m(\pa \PP) = \maK_{a}^m(\pa \PP; \vartheta) =
    \maK_{a}^m(\pa \PP; r_{\PP}) =
    r_{\PP}^{a-1/2} W^{m, 2}(\pa \PP).
\end{equation}
More precisely, let us notice that we can identify each edge with
$[0,1]$. Then $\maK_{a}^m(\pa \PP)$ consists of the functions $f :
\pa \PP \to \CC$ that, on each edge, are such that $t^k(1-t)^k
f^{(k)} \in L^2([0,1])$, $0 \le k \le m$ (here we identify that
edge with $[0, 1]$). This last condition is equivalent to
$[t(1-t)\pa_t]^k f \in L^2([0,1])$, $0 \le k \le m$.

\begin{proposition}\label{prop.poly}\
Let $\PP \subset \RR^2$ be a curvilinear polygonal domain and $P$
be a differential operator of order $m$ with coefficients in
$\CI(\overline{\PP})$. Then $P_\lambda := r_{\PP}^\lambda P
r_{\PP}^{-\lambda}$ defines a continuous family of bounded maps
$P_\lambda : \maK^{s}_{a}(\PP) \to \maK^{s-m}_{a-m}(\PP)$, for any
$s, a \in \RR$. Let $\PP'$ be $\overline{\PP}$ with the vertices
removed. Then $\CIc(\PP')$ is dense in $\maK^m_a(\PP)$. Also, the
restriction to the boundary extends to a continuous, surjective
trace map $\maK^s_a(\PP) \to \maK^{s-1/2}_{a-1/2}(\pa \PP)$. If $s
= 1$, then the kernel of the trace map is the closure of
$\CIc(\PP)$ in $\maK_{a}^1(\PP)$.
\end{proposition}

The above proposition, except maybe for the description of the
restrictions to the boundary, is well known in two dimensions. It
will serve as a model for the results in three dimensions that we
present in the last section.

\section{A regularity result\label{sec.5}}

We include in this section an application to the regularity of
boundary value problems, Theorem \ref{thm.reg1}. Its proof
is reduced to the Euclidean case using a partition of unity
argument and the tubular neighborhood theorem \ref{thm.tub.neigh},
both of which require some non-trivial input from differential
geometry.

Let us introduce some notation first that will be also useful in
the following. Let $\exp : TM_{0} \longrightarrow M_{0} \times
M_{0}$ be given by $\exp(v) := (x, \exp_{x}(v))$, $v \in T_xM_0$.
If $E$ is a real vector bundle with a metric, we shall denote by
$(E)_r$ the set of all vectors $v$ of $E$ with $|v|< r$. Let
$(M_0^2)_r := \{ (x, y), x, y \in M_0,\, d(x, y) < r\}$. Then the
exponential map defines a diffeomorphism $\exp :(TM_{0})_{r} \to
(M_0^2)_r$. We shall also need the admissible weight function
$\rho$ defined in Equation \eqref{eq.def.rho} and the weighted
Sobolev spaces $\rho^s \WkpO := \{\rho^s u, u \in \WkpO\}$
introduced in Equation \ref{eq.def.weights}.

Recall \cite{Taylor1}, Chapter 5, Equation (11.79), that a
differential operator $P$ of order $m$ is called {\em strongly
elliptic} if there exists $C > 0$ such that $Re \big(
\sigma^{(m)}(P)(\xi) \big) \ge C \|\xi\|^{m}$ for all $\xi$.

\begin{theorem}\label{thm.reg1}\
Let $\Omega \subset M$ be a regular open subset of the Lie manifold
$(M, \VV)$. Let $P \in \DiffV{2}(M)$ be an order $2$ strongly elliptic
operator on $M_0$ generated by $\VV$ and $s \in \RR$, $t \in \ZZ$, $1
< p < \infty$. Then there exists $C > 0$ such that, for any $u \in
\rho^{s}W^{1, p}(\Omega_0)$, $u\vert_{\pa \Omega_0} = 0$, we have
\begin{equation*}
    \|u\|_{\rho^{s}W^{t + 2, p}(\Omega_0)} \le
    C(\|Pu\|_{\rho^{s}W^{t, p}(\Omega_0)} +
    \|u\|_{\rho^{s}L^p(\Omega_0)}).
\end{equation*}
In particular, let $u \in \rho^{s}W^{1, p}(\Omega_0)$ be such that
$Pu \in \rho^{s}W^{t, p}(\Omega_0)$, and $u\vert_{\pa \Omega_0} =
0$, then $u \in \rho^{s}W^{t + 2, p}(\Omega_0)$.
\end{theorem}

\begin{proof}\
Note that, locally, this is a well known statement. In particular,
$\phi u \in W^{t + 2, p}(\Omega_0)$, for any $\phi \in \CIc(M_0)$.
The result will follow then if we prove that
\begin{equation}\label{eq.ell.est0}
    \|u\|_{\rho^{s}W^{t + 2, p}(M_0)} \le
    C(\|Pu\|_{\rho^{s}W^{t, p}(M_0)} +
    \|u\|_{\rho^{s}L^{p}(M_0)} )
\end{equation}
for any $u \in W^{t + 2, p}_{\mathrm{loc}}(\Omega_0)$. Here, of
course, $\|u\|_{\rho^{s}L^{p}(M_0)} = \|\rho^{-s}
u\|_{L^{p}(M_0)}$ (see Equation \eqref{eq.def.weights}).

Let $r < r_{\mathrm{inj}}(M_0)$ and let $\exp : (TM_0)_{r} \to
(M_0^2)_{r}$ be the exponential map. The statement is
trivially true for $t \le -2$, so we will assume $t \ge -1$ in
what follows. Also, we will assume first that $s = 0$. The general
case will be reduced to this one at the end. Assume first that
$\Omega_0 = M_0$.

Let $P_x$ be the differential operators defined on $B_{T_xM_0}(0,
r)$ obtained from $P$ by the local diffeomorphism $\exp :
B_{T_xM_0}(0, r) \to M_0$. We claim that there exists a constant $C >
0$, independent of $x \in M_0$ such that
\begin{equation}\label{eq.ell.reg}
    \|u\|_{W^{t+2, p}(T_xM_0)}^p \le C \big( \|P_xu\|_{W^{t,
    p}(T_xM_0)}^p + \|u\|_{L^{p}(T_xM_0)}^p \big),
\end{equation}
for any function $u \in \CIc(B_{T_xM_0}(0, r))$. This is seen as
follows. We can find a constant $C_x > 0$ with this property for
any $x \in M_0$ by the ellipticity of $P_x$. (For $p = 2$, a
complete proof can be found in \cite{Taylor1}, Propositions 11.10
and 11.16. For general $p$, the result can be proved as
\cite{Evans}, Theorem 1 in subsection 5.8.1, page~275.) Choose
$C_x$ to be the least such constant. Let $\pi : A \to M$ be the
extension of the tangent bundle of $M_0$, see Remark \ref{rem.A}
and let $A_x = \pi^{-1}(x)$. The family $P_x$, $x \in M_0$,
extends to a family $P_x$, $x \in M$, that is smooth in $x$. The
smoothness of the family $P_x$ in $x \in M$ shows that $C_x$ is
upper semi-continuous (\ie the set $\{C_x < \eta\}$ is
open for any $x$). Since $M$ is compact, $C_x$ will attain its
maximum, which therefore must be positive. Let $C$ be that
maximum value.

Let now $\phi_j$ be the partition of unity and $\psi_j$ be the
diffeomorphisms appearing in Equation \eqref{eq.nu.sp}, for some
$0 < \epsilon < r/6$. In particular, the partition of unity
$\phi_j$ satisfies the conditions of Lemma \ref{lemma.Shubin},
which implies that $\supp(\phi_j) \subset B(x_j, 2\epsilon)$ and
the sets $B(x_j, 4\epsilon)$ form a covering of $M_0$ of finite
multiplicity.  Let $\eta_j = 1$ on the support of $\phi_j$,
$\supp(\eta_j) \subset B(x_j, 4\epsilon)$. We then have
\begin{multline*}
    \nu_{t + 2, p}(u)^p :=  \sum_j \|(\phi_j u) \circ
    \psi_j^{-1}\|_{W^{t+2,p}(\RR^n)}^p \\ \le C \sum_j \Big( \|P_x
    (\phi_j u)\|_{W^{t, p}(T_xM_0)}^p + \|\phi_j u\|_{L^{
    p}(T_xM_0)}^p \Big) \\ \le C \sum_j \Big( \| \phi_j P_x u\|_{W^{t,
    p}(T_xM_0)}^p + \|[P_x, \phi_j] u \|_{W^{t, p}(T_xM_0)}^p +
    \|\phi_j u\|_{L^{p}(T_xM_0)}^p \Big) \\ \le C \sum_j \Big( \|
    \phi_j P_x u\|_{W^{t, p}(T_xM_0)}^p + \| \eta_j u \|_{W^{t+1,
    p}(T_xM_0)}^p + \|\phi_j u\|_{L^{p}(T_xM_0)}^p \Big) \\ \le C
    \big( \nu_{t, p}(P u)^p + \nu_{t+1}(u)^p \big).
\end{multline*}
The equivalence of the norm $\nu_{s, p}$ with the standard norm on
$W^{s, p}(M_0)$ (Propositions \ref{prop.alt.desc} and
\ref{prop.alt.desc1}) shows that $\|u\|_{W^{t + 2, p}(M_0)} \le C(
\|Pu\|_{W^{t, p}(M_0)} + \|u\|_{W^{t + 1, p}(M_0)} )$, for any $t
\ge -1$. This is known to imply
\begin{equation}\label{eq.ell.est}\
    \|u\|_{W^{t + 2, p}(M_0)} \le C( \|Pu\|_{W^{t, p}(M_0)} +
    \|u\|_{L^{p}(M_0)} )
\end{equation}
by a boot-strap procedure, for any $t \ge -1$. This proves our
statement if $s = 0$ and $\Omega_0 = M_0$.

The case of arbitrary domains $\Omega_0$ follows in exactly the
same way, but using a product type metric in a neighborhood of
$\patop \Omega_0$ and the analogue of Equation \eqref{eq.ell.reg}
for a half-space, which shows that Equation \eqref{eq.ell.est0}
continues to hold for $M_0$ replaced with $\Omega_0$.

The case of arbitrary $s\in \RR$ is obtained by applying Equation
\eqref{eq.ell.est} to the elliptic operator $\rho^{-s} P \rho^{s}
\in \DiffV{2}(M)$ and to the function $\rho^{-s}u \in W^{k,
p}(\Omega_0)$, which then gives Equation \eqref{eq.ell.est0} right
away.
\end{proof}

For $p = 2$, by combining the above theorem with Theorem
\ref{theorem.res}, we obtain the following corollary.

\begin{corollary}\label{cor.reg1}\
We keep the assumptions of Theorem \ref{thm.reg1}. Let $u \in
\rho^{s}H^{1}(\Omega_0)$ be such that $Pu \in
\rho^{s}H^{t}(\Omega_0)$ and $u \vert_{\pa \Omega_0} \in
\rho^{s}H^{t + 3/2}(\Omega_0)$, $s \in \RR$, $t \in \ZZ$. Then $u
\in \rho^{s}H^{t + 2}(\Omega_0)$ and
\begin{equation}\label{eq.ell.estB}
    \|u\|_{\rho^{s}H^{t + 2}(\Omega_0)} \le
    C(\|Pu\|_{\rho^{s}H^{t}(\Omega_0)} +
    \|u\|_{\rho^{s}L^2(\Omega_0)} +
    \|u\vert_{\pa \Omega_0}\|_{\rho^{s}H^{t + 3/2}(\Omega_0)}).
\end{equation}
\end{corollary}

\begin{proof}\
For $u\vert_{\pa \Omega_0} = 0$, the result follows from Theorem
\ref{thm.reg1}. In general, choose a suitable $v \in
H^{t+2}(\Omega_0)$ such that $v\vert_{\pa \Omega_0} = u\vert_{\pa
\Omega_0}$, which is possible by Theorem \ref{theorem.res}. Then
we use our result for $u-v$.
\end{proof}

\section{Polyhedral domains in three dimensions\label{sec.PD}}

We now include an application of our results to polyhedral domains
$\PP \subset \RR^3$. A~{\em polyhedral domain in} $\PP \subset
\RR^3$ is a bounded, connected open set such that $\pa \PP = \pa
\overline{\PP}= \bigcup \overline{D}_j$ where
\begin{itemize}
\item each $D_j$ is a polygonal domain with straight edges contained
in an affine $2$-dimensional subspace of $\RR^3$
\item each edge is contained in exactly two closures of polygonal
domains $\overline{D}_j$.
\end{itemize}
(See Subsection \ref{ssec.poly} for the definition of a polygonal
domain.)

The vertices of the polygonal domains $D_j$ will form the vertices of
$\PP$. The edges of the polygonal domains $D_j$ will form the edges of
$\PP$.  For each vertex $P$ of $\PP$, we choose a small open ball
$V_P$ centered in $P$. We assume that the neighborhoods $V_P$ are
chosen to be disjoint.  For each vertex $P$, there exists a unique
closed polyhedral cone $C_P$ with vertex at $P$, such that
$\overline{\PP} \cap V_P = C_P \cap V_P$. Then $\PP\subset\bigcup
C_P$.

We now proceed to define canonical weight functions of
$\PP$ in analogy with the definition of canonical weights of
curvilinear polygonal domains, Definition \ref{def.can.weight}. We
want to define first a continuous function $r_{\PP} :
\overline{\Omega} \to [0, \infty)$ that is positive and
differentiable outside the edges. Let $\vartheta(x)$ be the
distance from $x$ to the edges of $\PP$, as before. We want
$r_{\PP}(x) = \vartheta(x)$ close to the edges but far from the
vertices and we want the quotients $r_{\PP}(x)/\vartheta(x)$ and
$\vartheta(x)/r_{\PP}(x)$ to extend to continuous functions on
$\overline{\Omega}$. Using a smooth partition of unity, in order
to define $r_{\PP}$, we need to define it close to the vertices.

Let us then denote by $\{P_k\}$ the set of vertices of
$\PP$. Choose a continuous function $r: \overline{\PP} \to [0,
\infty)$ such that $r(x)$ is the distance from $x$ to the vertex
$P$ if $x \in V_P \cap \PP$, and such that $r(x)$ is
differentiable and positive on $\overline{\PP} \smallsetminus
\{P_k\}$.  Let $S^2$ be the unit sphere centered at $P$ and let
$r_P$ be a canonical weight associated to the curvilinear polygon
$C_P \cap S^2$ (see Definition \ref{def.can.weight}). We extend
this function to $C_P$ to be constant along the rays, except at
$P$, where $r_P(P) = 0$. Finally, we let $r_\PP(x) = r(x) r_P(x)$,
for $x$ close to $P$. Then a {\em canonical weight of $\PP$} is
any function of the form $\psi r_{\PP}$, where $\psi$ is a smooth,
nowhere vanishing function on $\overline{\PP}$.

For any canonical weight $r_{\PP}$, we then we have the following
analogue of Equation \eqref{eq.same.Sobo}
\begin{equation}\label{eq.same.Sobo2}
    \maK_a^m(\PP) := \maK_a^m(\PP; \vartheta)
    = \maK_a^m(\PP; r_{\PP}).
\end{equation}

Let us define, for every vertex $P$ of $\PP$, a spherical
coordinate map $\Theta_P : \PP \smallsetminus \{P\} \to S^2$ by
$\Theta_P(x) = |x - P|^{-1}(x - P)$. Then, for each edge $e =
[AB]$ of $\PP$ joining the vertices $A$ and $B$, we define a
generalized cylindrical coordinate system $(r_e, \theta_e, z_e)$
to satisfy the following properties:
\begin{enumerate}[(i)]
\item\ $r_e(x)$ be the distance from $x$ to the line containing $e$.
\item\ $A$ as the origin (\ie $r_e(A) = z_e(A) = 0)$,
\item\ $\theta_e = 0$ on one of the two faces containing $e$, and
\item\ $z_e \ge 0$ on the edge $e$.
\end{enumerate}
Let $\psi : S^2 \to [0, 1]$ be a smooth function on the unit sphere
that is equal to $1$ in a neighborhood of $(0, 0, 1) = \{\phi = 0\}
\cap S^2$ and is equal to $0$ in a neighborhood of $(0, 0, -1) =
\{\phi = \pi\} \cap S^2$. Then we let
\begin{equation*}
    \tilde \theta_e(x) = \theta_e(x) \psi(\Theta_A(x))
    \psi(-\Theta_B(x))
\end{equation*}
where $\theta_e(x)$ is the $\theta$ coordinate of $x$ in a
cylindrical coordinate system $(r, \theta, z)$ in which the point $A$
corresponds to the origin (\ie $r = 0$ and $z = 0$) and the edge $AB$
points in the positive direction of the $z$ axis (\ie $B$ corresponds
to $r = 0$ and $z > 0$). By choosing $\psi$ to have support small
enough in $S^2$ we may assume that the function $\tilde \theta_e$ is
defined everywhere on $\PP \smallsetminus e$. (This is why we need the
cut-off function $\psi$.)

We then consider the function
\begin{equation*}
    \Phi : \PP \to \RR^N, \quad \Phi(x) = (x, \Theta_P(x), r_e(x),
    \tilde \theta_e(x)),
\end{equation*}
with $N = 3 + 3n_v + 2n_e$, $n_v$ being the number of vertices of
$\PP$ and $n_e$ being the number of edges of $\PP$. Finally, we
define $\Sigma(\PP)$ to be the closure of $\Phi(\PP)$ in $\RR^N$.
Then $\Sigma(\PP)$ is a manifold with corners that can be endowed
with the structure of a Lie manifold with true boundary as
follows. (Recall that a Lie manifold with boundary $\Sigma$ is the
closure $\overline{\Omega}$ of a regular open subset $\Omega$ in a
Lie manifold $M$ and the {\em true boundary} of $\Sigma$ is the
topological boundary $\patop \Omega$.) The true boundary $\patop
\Sigma(\Omega)$ of $\Sigma(\Omega)$ is defined as the union of the
closures of the faces $D_j$ of $\PP$ in $\Sigma(\PP)$. (Note that
the closures of $D_j$ in $\Sigma(\PP)$ are disjoint.) We can then
take $M$ to be the union of two copies of $\Sigma(\PP)$ with the
true boundaries identified (\ie the double of $\Sigma(\PP)$) and
$\Omega = \Sigma(\PP) \smallsetminus \patop \Sigma(\PP)$. In
particular, $\Omega_0 := \Omega \cap M_0$ identifies with $\PP$.

To complete the definition of the Lie manifold with true boundary
on $\Sigma(\PP)$, we now define the structural Lie algebra of
vector fields $\VV(\PP)$ of $\Sigma(\PP)$ by
\begin{equation}
    \VV(\PP) := \{ r_\PP(\phi_1 \pa_1 + \phi_2 \pa_2 + \phi_3
    \pa_3), \, \phi_j \in \CI(\Sigma(\PP)) \}.
\end{equation}
(Here $\pa_j$ are the standard unit vector fields. Also, the
vector fields in $\VV(\PP)$ are determined by their restrictions
to $\PP$.) This is consistent with the fact that $\patop
\Sigma(\PP)$, the true boundary of $\Sigma(\PP)$, is defined as
the union of the boundary hyperfaces of $\Sigma(\PP)$ to which not
all vector fields are tangent. This completes the definition of
the structure of Lie manifold with boundary on $\Sigma(\PP)$.

The function $r_{\PP}$ is easily seen to be an admissible weight
on $\Sigma(\PP)$. It hence satisfies
\begin{equation*}
    r_{\PP} (\pa_j r_{\PP}) = r_{\PP}
    \frac{\pa r_{\PP}}{\pa x_j} \in \CI(\Sigma(\PP)),
\end{equation*}
which is equivalent to the fact that $\VV(\PP)$ is a Lie algebra.
This is the analogue of Equation \eqref{eq.vf2}.

To check that $\Sigma(\PP)$ is a Lie manifold, let us notice first
that $g = r_{\PP}^{-2} g_E$ is a compatible metric on
$\Sigma(\PP)$, where $g_E$ is the Euclidean metric on $\PP$. Then,
let us denote by $\nu$ the outer unit normal to $\PP$ (where it is
defined), then $r_{\PP} \pa_\nu$ is the restriction to $\patop
\Sigma(\Omega)$ of a vector field in $\VV(\PP)$. Moreover $r_{\PP}
\pa_\nu$ is of length one and orthogonal to the true boundary in
the compatible metric $g = r_{\PP}^{-2} g_E$.

The definition of $\VV(\PP)$ together with our definition of
Sobolev spaces on Lie manifolds using vector fields shows that
\begin{equation}
    \maK_a^m(\PP) = r_{\PP}^{a-3/2} W^{m, 2}(\Sigma(\PP))
    = r_{\PP}^{a-3/2}H^m(\Sigma(\PP)).
\end{equation}

The induced Lie manifold structure on $\Sigma(\PP)$ consists of
the vector fields on the faces $D_j$ that vanish on the boundary
of $D_j$. The Soblev spaces on the boundary are
\begin{equation}
    \maK_a^m(\pa \PP) = r_{\PP}^{a-1}
    W^{m, 2}(\patop \Sigma(\PP)) = r_{\PP}^{a-1}H^m(\patop
    \Sigma(\PP)).
\end{equation}
The factors $-3/2$ and $-1$ in the powers of $r_{\PP}$ appearing
in the above two equations are due to the fact that the volume
elements on $\PP$ and $\Sigma(\PP)$ differ by these factors.

If $P$ is an order $m$ differential operator with smooth
coefficients on $\RR^3$ and $\PP \subset \RR^3$ is a polyhedral
domain, then $r_{\PP}^m P \in \DiffV{m}(\Sigma(\PP))$, by Equation
\eqref{eq.vf1}. However, in general, $r_{\PP}^m P$ will not define
a smooth differential operator on $\overline{\PP}$.

In particular, we have the following theorem, which is a direct
analog of Proposition \ref{prop.poly}, if we replace ``vertices''
with ``edges:''

\begin{theorem}\label{thm.3D}\
Let $\PP \subset \RR^3$ be a polyhedral domain and $P$ be a
differential operator of order $m$ with coefficients in
$\CI(\overline{\PP})$. Then $P_\lambda := r_{\PP}^\lambda P
r_{\PP}^{-\lambda}$ defines a continuous family of bounded maps
$P_\lambda : \maK^{s}_{a}(\PP) \to \maK^{s-m}_{a-m}(\PP)$, for any
$s, a \in \RR$. Let $\PP'$ be $\overline{\PP}$ with the {\em
edges} removed. Then $\CIc(\PP')$ is dense in $\maK^m_a(\PP)$.
Also, the restriction to the boundary extends to a continuous,
surjective trace map $\maK^s_a(\PP) \to \maK^{s-1/2}_{a-1/2}(\pa
\PP)$. If $s = 1$, then the kernel of the trace map is the closure
of $\CIc(\PP)$ in $\Kond{1}{a}(\PP)$.
\end{theorem}

See \cite{BNZ-I} for applications of these results,
especially of the above theorem.

Theorem \ref{thm.reg1} and the results of this section
immediately lead to the proof of Theorem \ref{theorem.3Dreg}
formulated in the Introduction.

\section{A non-standard boundary value problem\label{sec.NS}}

We present in this section a non-standard boundary value problem on
a smooth manifold with boundary. Let $\maO$ be a smooth manifold
with boundary. We shall assume that $\maO$ is connected and that
the boundary is not empty.

Let $r: \overline{\maO} \to [0, \infty)$ be a smooth function that
close to the boundary is equal to the distance to the boundary and
is $>0$ on $\maO$. Then we recall \cite{Donnelly} that there
exists a constant depending only on $\maO$ such that
\begin{equation}\label{eq.Hardy}
    \int_{\maO} r^{-2}|u(x)|^2 dx \le C \int_{\maO}
    |\nabla u(x)|^2 dx
\end{equation}
for any $u \in H^1(\maO)$ that vanishes at the boundary. If we
denote, as in Equation \eqref{def.t.c.a0},
\begin{equation*}
    \Kond{m}{a}(\maO; r) := \{ u \in L^2_{loc}(\maO), \
    r^{|\alpha|-a} \pa^{\alpha} u \in L^2(\maO), \ \
    |\alpha| \le m \}, \quad m \in \NN\cup\{0\},\, a \in \RR,
\end{equation*}
with norm $\|\, \cdot\, \|_{\Kond{m}{a}}$, the Equation
\eqref{eq.Hardy} implies that $\|u\|_{\Kond{1}{1}} \le C\|\nabla
u\|_{L^2}$.

Let $M = \overline{\maO}$ with the structural Lie algebra of
vector fields
\begin{equation*}
    \VV =\VV_0 := \{X, X = 0 \text{ at } \pa \maO\} = r\Gamma(M; TM),
\end{equation*}
(see Example \ref{ex1}). Recall from Subsection \ref{ssec.Diff} that
$\DiffV{m}(M)$ is the space of order $m$ differential operators on $M$
generated by multiplication with functions in $\CI(M)$ and by
differentiation with vector fields $X \in \VV$.  It follows that
\begin{equation}\label{eq.One}
    r^m P \in \DiffV{m}(M)
\end{equation}
for any differential operator $P$ of order $m$ with smooth
coefficients on $M$.

\begin{lemma}\label{lemma.One}\
The pair $(M, \VV)$ is a Lie manifold with $M_0 = \maO$ satisfying
\begin{equation}
    \Kond{m}{a}(\maO; r) = r^{a-n/2}H^m(M).
\end{equation}
If $P$ is a differential operator with smooth coefficients on $M$,
then $r^m P$ is a differential operator generated by $\VV$, and
hence $P_\lambda := r^\lambda P r^{-\lambda}$ gives rise to a
continuous family of bounded maps $P_{\lambda} : \Kond{s}{a}(\maO;
r) \to \Kond{s-m}{a-m}(\maO; r)$.
\end{lemma}

Because of the above lemma, it makes sense to define
$\Kond{s}{a}(\maO; r) = r^{a-n/2}H^s(M)$, for all $s, a \in \RR$, with
norm denoted $\|\cdot\|_{\Kond{s}{a}}$. The regularity result (Theorem
\ref{thm.reg1}) then gives

\begin{lemma}\label{lemma.Two}\
Let $P$ be an order $m$ elliptic differential operator with smooth
coefficients defined in a neighborhood of $M = \overline{\maO}$.
Then, for any $s, t \in \RR$, there exists $C = C_{st} > 0$ such that
\begin{equation*}
    \|u\|_{\Kond{s}{a}} \le
    C\big(\|Pu\|_{\Kond{s-m}{a-m}} +
    \|u\|_{\Kond{t}{a}}\big).
\end{equation*}
In particular, let $u \in \Kond{t}{a}(\maO; r)$ be such that $Pu \in
\Kond{s-m}{a-m}(\maO; r)$, then $u \in \Kond{s}{a}(\maO; r)$.  The
same result holds for elliptic systems.
\end{lemma}

\begin{proof}\
We first notice that $r^mP \in \DiffV{m}(M)$ is an {\em elliptic
operator} in the usual sense (that is, its principal symbol
$\sigma^{(m)}(r^mP)$ does not vanish outside the zero section of
$A^*$). For this we use that $\sigma^{(m)}(r^mP) = r^m
\sigma^{(m)}(P)$ and that $A^*$ is defined such that multiplication by
$r^m$ defines an isomorphism $\CI(T^*M) \to \CI(A^*)$ that maps order
$m$ elliptic symbols to elliptic symbols. Then the proof is exactly
the same as that of Theorem \ref{thm.reg1}, except that we do not need
strong ellipticity, because we do not have boundary conditions (and
hence we have no condition of the form $u = 0$ on the boundary).
\end{proof}

An alternative proof of our lemma is obtained using
pseudodifferential operators generated by $\VV$ \cite{aln2} and
their $L^p$--continuity.

\begin{theorem}\
There exists $\eta > 0$ such that $\Delta : \Kond{a+1}{s+1}(\maO;
r) \to \Kond{a-1}{s-1}(\maO; r)$ is an isomorphism for all $s \in
\RR$ and all $|a| < \eta$.
\end{theorem}

\begin{proof}\
The proof is similar to that of Theorem 2.1 in \cite{BNZ1}, so we
will be brief. Consider
\begin{equation*}
    B: \Kond{1}{1}(\maO; r) \times \Kond{1}{1}(\maO; r) \to \CC,
    \quad B(u, v) = \int_\maO \nabla u \cdot \nabla \overline{v} dx.
\end{equation*}
Then $|B(u, v)| \le \|u\|_{\Kond{1}{1}}\|v\|_{\Kond{1}{1}}$, so
$B$ is continuous.

On the other hand, by Equation \eqref{eq.Hardy}, $B(u, u) \ge
\theta \|u\|_{\Kond{1}{1}}^2$,
for all $u$ with compact support on
$\maO$ and for some $\theta > 0$ independent of $u$. Since
$\CIc(\maO)$ is dense in $\Kond{1}{1}(\maO; r)$, by Theorem
\ref{theorem.dense}, the Lax-Milgram Lemma can be used to conclude
that
\begin{equation*}
    \Delta : \Kond{1}{1}(\maO; r) \to \Kond{-1}{-1}(\maO; r)
    := \Kond{1}{1}(\maO; r)^*
\end{equation*}
is an isomorphism. Since multiplication by $r^a :
\Kond{1}{1}(\maO; r) \to \Kond{1}{a+1}(\maO; r)$ is an isomorphism
and the family $r^{a}\Delta r^{-a}$ depends continuously on $a$ by
Lemma \ref{lemma.One}, we obtain that $\Delta :
\Kond{1}{a+1}(\maO; r) \to \Kond{-1}{a-1}(\maO; r)$ is an
isomorphism for $|a|<\eta$, for some $\eta>0$ small enough.

Fix now $a$, $|a| < \eta$. We obtain that $\Delta :
\Kond{s+1}{a+1}(\maO; r) \to \Kond{s-1}{a-1}(\maO; r)$ is a
continuous, injective map, for all $s \ge 0$. The first part of the
proof (for $a = 0$) together with the regularity result of Lemma
\ref{lemma.Two} show that this map is also surjective. The Open
Mapping Theorem therefore completes the proof for $s \ge 0$. For $s
\le 0$, the result follows by considering duals.
\end{proof}

It can be shown as in \cite{BNZ1} that $\eta$ is the least value
for which $\Delta : \Kond{1}{\eta+1}(\maO; r) \to
\Kond{-1}{\eta-1}(\maO; r)$ is not Fredholm. This, in principle,
can be decided by using the Fredholm conditions in
\cite{NistorPrden} that involve looking at the $L^2$ invertibility
of the same differential operators when $M$ is the half-space
$\{x_{n+1} \ge 0\}$. See also \cite{Giroire} for some non-standard
boundary value problems on exterior domains in weighted Sobolev
spaces.

\section{Pseudodifferential operators\label{sec.PO}}

We now recall the definition of pseudodifferential operators on
$M_0$ generated by a Lie structure at infinity $(M, \VV)$ on
$M_0$.

\subsection{Definition}

We fix in what follows a compatible Riemannian metric $g$ on $M_0$
(that is, a metric coming by restriction from a metric on the
bundle $A \to M$ extending $TM_0$), see Section \ref{sec.LI}. In
order to simplify our discussion below, we shall use the metric
$g$ to trivialize all density bundles on $M$. Recall that $M_0$
with the induced metric is complete \cite{aln1}. Also, recall that
$A \to M$ is a vector bundle such that $\VV = \Gamma(A)$.

Let $\exp_x : T_x M_0 \to M_0$ be the exponential map, which is
everywhere defined because $M_0$ is complete. We let
\begin{equation}\label{eq.Phi}
    \Phi : TM_{0} \longrightarrow M_{0} \times M_{0}, \quad
    \Phi(v) := (x, \exp_{x}(-v)), \; v \in T_xM_0,
\end{equation}

If $E$ is a real vector bundle with a metric, we shall denote by
$(E)_r$ the set of all vectors $v$ of $E$ with $|v|< r$. Let
$(M_0^2)_r := \{ (x, y), x, y \in M_0,\, d(x, y) < r\}$. Then the
map $\Phi$ of Equation \eqref{eq.Phi} restricts to a
diffeomorphism $\Phi :(TM_{0})_{r} \to (M_0^2)_r$, for any $0 < r
< r_{\mathrm{inj}}(M_0)$, where $r_{\mathrm{inj}}(M_0)$ is the
injectivity radius of $M_0$, which was assumed to be positive. The
inverse of $\Phi$ is of the form
\begin{equation*}
    (M_0^2)_{r} \ni (x,y) \longmapsto
    (x,\tau(x,y))\in (TM_{0})_{r}\,.
\end{equation*}

We shall denote by $S^m_{1,0}(E)$ the space of symbols of order
$m$ and type $(1,0)$ on $E$ (in H\"ormander's sense) and by
$S^m_{cl}(E)$ the space of classical symbols of order $m$ on $E$
\cite{hoe67, NS, Taylor, Taylor2}. See \cite{aln2} for a review of
these spaces of symbols in our framework.

Let $\chi \in \CI(A^*)$ be a smooth function that is equal to $1$
on $(A^*)_{r}$ and is equal to $0$ outside $(A^*)_{2r}$, for some
$r < r_{\mathrm{inj}}(M_0)/3$. Then, following \cite{aln2}, we
define
\begin{equation*}
    q(a)u(x) = (2\pi )^{-n} \int_{T^{*}M_{0}} e^{i\tau(x,y)\cdot\eta}
    \chi(x,\tau(x,y)) a(x,\eta)u(y)\, d \eta \,dy \,.
\end{equation*}
This integral is an oscillatory integral with respect to the
symplectic measure on $T^*M_0$ \cite{hor3}. Alternatively, we
consider the measures on $M_0$ and on $T^*_xM_0$ defined by some
choice of a metric on $A$ and we integrate first in the fibers
$T^*_xM_0$ and then on $M_0$. The map $\sigma_{tot} :
S_{1,0}^{m}(A^{*}) \to \Psi^{m}(M_{0})/\Psi^{-\infty}(M_{0})$,
\begin{equation*}
    \sigma_{tot}(a) := q(a) + \Psi^{-\infty}(M_{0})
\end{equation*}
is independent of the choice of the function $\chi \in
\CIc((A)_{r})$ \cite{aln2}.

We now enlarge the class of order $-\infty$ operators that we
consider. Any $X \in \VV = \Gamma(A)$ generates a global flow
$\Psi_{X} : \RR\times M\rightarrow M$ because $X$ is tangent to
all boundary faces of $M$ and $M$ is compact.  Evaluation at $t =
1$ yields a diffeomorphism
\begin{equation}\label{eq.def.psiX}
    \psi_X := \Psi_{X}(1,\cdot) : M \rightarrow M.
\end{equation}

We now define the pseudodifferential calculus on $M_0$ that we will
consider following \cite{aln2}. See \cite{lmn1, lmn2, Monthubert, NWX}
for the connections between this calculus and groupoids.

\begin{definition}\label{def.psi.MA}\
Fix $0< r < r_{\mathrm{inj}}(M_0)$ and $\chi \in \CIc((A)_{r})$
such that $\chi = 1$ in a neighborhood of $M \subseteq A$.  For $m
\in \RR$, the space $\Psi_{1,0,\VV}^{m}(M_{0})$ of {\em
pseudodifferential operators generated by the Lie structure at
infinity $(M, \VV)$} is defined to be the linear space of
operators $\CIc(M_0) \rightarrow \CIc(M_{0})$ generated by $q(a)$,
$a \in S_{1,0}^m(A^*)$, and $q(b)\psi_{X_1}\ldots \psi_{X_k}$, $b
\in S^{-\infty}(A^*)$ and $X_j \in \Gamma(A)$, $\forall j$.

Similarly, the space $\Psi_{cl,\VV}^{m}(M_{0})$ of {\em classical
pseudodifferential operators generated by the Lie structure at
infinity $(M, \VV)$} is obtained by using classical symbols $a$ in
the construction above.
\end{definition}

We have that $\Psi_{cl,\VV}^{-\infty}(M_{0}) =
\Psi_{1,0,\VV}^{-\infty}(M_{0}) =: \Psi_{\VV}^{-\infty}(M_{0})$
(we dropped some subscripts).

\subsection{Properties}

 We now review some properties of the
operators in $\Psi^{m}_{1, 0, \VV}(M_{0})$ and $\Psi^{m}_{cl,
\VV}(M_{0})$ from \cite{aln2}. These properties will be used
below. Let $\Psi^{\infty}_{1, 0, \VV}(M_{0}) = \bigcup_{m \in \ZZ}
\Psi^{m}_{1, 0, \VV}(M_{0})$ and $\Psi^{\infty}_{cl, \VV}(M_{0}) =
\bigcup_{m \in \ZZ} \Psi^{m}_{cl, \VV}(M_{0})$.

First of all, each operator $P \in \Psi^{m}_{1, 0, \VV}(M_{0})$
defines continuous maps $\CIc(M_{0}) \rightarrow \CI(M_{0})$, and
$\CI(M) \rightarrow \CI(M)$, still denoted by $P$. An operator $P
\in \Psi^{m}_{1, 0, \VV}(M_{0})$ has a distribution kernel $k_P$
in the space $I^m(M_0 \times M_0, M_0)$ of distributions on $M_0
\times M_0$ that are conormal of order $m$ to the diagonal, by
\cite{hor3}. If $P = q(a)$, then $k_P$ has support in $(M_0 \times
M_0)_r$. If we extend the exponential map $(TM_0)_r \to M_0 \times
M_0$ to a map $A \to M$, then the distribution kernel of $P =
q(a)$ is the restriction of a distribution, also denoted $k_P$ in
$I^m(A, M)$.

If $\maP$ denotes the space of polynomial symbols on $A^*$ and
$\Diff{M_0}$ denotes the algebra of differential operators on
$M_0$, then
\begin{equation}
    \Psi^{\infty}_{1, 0, \VV}(M_{0}) \cap \Diff{M_0}
    = \DiffV{\infty}{(M)} = q(\maP).
\end{equation}

The spaces $\Psi_{1, 0, \VV}^m(M_{0})$ and $\Psi_{1, 0,
\VV}^m(M_{0})$ are independent of the choice of the metric on $A$
and the function $\chi$ used to define it, but depend, in general,
on the Lie structure at infinity $(M,A)$ on $M_0$. They are also
closed under multiplication, which is a quite non-trivial fact.

\begin{theorem}\label{thm.alg}\
The spaces $\Psi_{1,0, \VV}^{\infty}(M_{0})$ and $\Psi_{cl,
\VV}^{\infty}(M_{0})$ are filtered algebras that are closed under
adjoints.
\end{theorem}

For $\Psi_{1,0, \VV}^m(M_{0})$, the meaning of the above theorem
is that
\begin{equation*}
    \Psi_{1,0, \VV}^m(M_{0}) \Psi_{1, 0, \VV}^{m'}(M_{0})
    \subseteq \Psi_{1, 0, \VV}^{m + m'}(M_{0}) \text{ and }
    \big(\Psi_{1,0, \VV}^m(M_{0})\big)^* = \Psi_{1,0,
    \VV}^m(M_{0})
\end{equation*}
for all $m,m' \in \CC \cup \{-\infty\}.$

The usual properties of the principal symbol remain true.

\begin{proposition}\label{prop.princ.symb}\
The principal symbol establishes isomorphisms
\begin{equation}\label{eq.def.princ}
    \sigma^{(m)} : \Psi_{1,0,\VV}^{m}(M_{0}) /
    \Psi_{1,0,\VV}^{m-1}(M_{0}) \to S^m_{1,0}(A^*) /
    S^{m-1}_{1,0}(A^*)
\end{equation}
and
\begin{equation}\label{eq.def.princ'}
    \sigma^{(m)} : \Psi_{cl,\VV}^{m}(M_{0}) /
    \Psi_{cl,\VV}^{m-1}(M_{0}) \to S^m_{cl}(A^*) /
    S^{m-1}_{cl}(A^*).
\end{equation}
Moreover, $\sigma^{(m)}(q(a)) = a + S^{m-1}_{1, 0}(A^*)$ for any
$a \in S^{m}_{1, 0}(A^*)$ and $\sigma^{(m+m')}(PQ) =
\sigma^{(m)}(P)\sigma^{(m')}(Q)$, for any $P \in
\Psi_{1,0,\VV}^{m}(M_{0})$ and $Q \in \Psi_{1,0,\VV}^{m'}(M_{0})$.
\end{proposition}

We shall need also the following result.

\begin{proposition} \label{prop.auto}\
Let $\rho$ be a defining function of some hyperface of $M$. Then
$\rho^{s} \Psi_{1,0, \VV}^m(M_{0}) \rho^{-s} = \Psi_{1, 0,
\VV}^m(M_{0})$ and $\rho^{s} \Psi_{cl, \VV}^m(M_{0}) \rho^{-s} =
\Psi_{cl, \VV}^m(M_{0})$ for any $s \in \CC$.
\end{proposition}

\subsection{Continuity on $W^{s,p}(M_0)$}

The preparations above will allow us to prove the continuity of
the operators $P \in \Psi_{1, 0, \VV}^m(M_{0})$ between suitable
Sobolev spaces. This is the main result of this section.  Some of
the ideas and constructions in the proof below have already been
used in \ref{thm.reg1}, which the reader may find convenient to
review first. Let us recall from Equation
\eqref{eq.def.rho} that an admissible weight $\rho$ of $M$ is a
function of the form $\rho := \prod_H \rho_H ^{a_H}$, where $a_H
\in \RR$ and $\rho_H$ is a defining function of $H$.

\begin{theorem}\label{thm.Wkp.cont}\
Let $\rho$ be an admissible weight of $M$ and let $P \in
\Psi_{1,0, \VV}^m(M_{0})$ and $p \in (0,\infty)$. Then $P$ maps
$\rho^rW^{s, p}(M_0)$ continuously to $\rho^rW^{s-m, p}(M_0)$ for
any $r, s \in \RR$.
\end{theorem}

\begin{proof}\
We have that $P$ maps $\rho^rW^{s, p}(M_0)$ continuously to
$\rho^rW^{s-m, p}(M_0)$ if, and only if, $\rho^{-r}P\rho^r$ maps
$W^{s, p}(M_0)$ continuously to $W^{s-m, p}(M_0)$. 
By Proposition~\ref{prop.auto} it is therefore enough to check 
our result for $r=0$.

We shall first prove our result if the Schwartz kernel of $P$ has
support close enough to the diagonal. To this end, let us choose
$\epsilon < r_{\mathrm{inj}}(M_0)/9$ and assume that the
distribution kernel of $P$ is supported in the set
$(M_0^2)_{\epsilon} := \{(x, y), d(x, y) < \epsilon\} \subset
M_0^2$. This is possible by choosing the function $\chi$ used to
define the spaces $\Psi_{1, 0, \VV}^m(M_0)$ to have support in the
set $(M_0^2)_{\epsilon}$. There will be no loss of generality then
to assume that $P = q(a)$.

Then choose a smooth function $\eta : [0, \infty) \to [0, 1]$,
$\eta(t) = 1$ if $t \le 6\epsilon$, $\eta(t) = 0$ if $t \ge
7\epsilon$. Let $\psi_x : B(x, 8\epsilon) \to B_{T_xM_0}(0,
8\epsilon)$ denote the normal system of coordinates induced by the
exponential maps $\exp_x : T_x M_0 \to M_0$. Denote $\pi : A \to
M$ be the natural (vector bundle) projection and
\begin{equation}
    B := A \times_{M}A := \{(\xi_1, \xi_2) \in A \times A, \pi(\xi_1)
    = \pi(\xi_2) \},
\end{equation}
which defines a vector bundle $B \to M$. In the language of vector
bundles, $B := A \oplus A$. For any $x \in M_0$, let $\eta_x$
denote the function $\eta \circ \exp_x$, and consider the operator
$\eta_x P \eta_x$ on $B(x, 13\epsilon)$. The diffeomorphism
$\psi_x$ then will map this operator to an operator $P_x$ on
$B_{T_xM_0}(0, 8\epsilon)$. Then $P_x$ maps continuously $W^{s,
p}(T_xM_0) \to W^{s-m, p}(T_xM_0)$, by the continuity of
pseudodifferential operators on $\RR^n$ \cite[XIII, \S 5]{Taylor3}
or \cite{Stein}.

The distribution kernel $k_x$ of $P_x$ is a distribution with
compact support on
\begin{equation*}
    T_xM_0 \times T_xM_0 = A_x \times A_x = B_x
\end{equation*}
If $P = q(a) \in \Psi_{1, 0, \VV}^m(M_0)$, then the distributions
$k_x$ can be determined in terms of the distribution $k_P \in
I^m(A, M)$ associated to $P$. This shows that the distributions
$k_x$ extend to a smooth family of distributions on the fibers of
$B \to M$. From this, it follows that the family of operators $P_x
: W^{s, p}(A_x) \to W^{s - m, p}(A_x)$, $x \in M_0$, extends to a
family of operators defined for $x \in M$ (recall that $A_x =
T_xM_0$ if $x \in M_0$). This extension is obtained by extending
the distribution kernels. In particular, the resulting family
$P_x$ will depend smoothly on $x \in M$. Since $M$ is compact, we
obtain, in particular, that the norms of the operators $P_x$ are
uniformly bounded for $x \in M_0$.

By abuse of notation, we shall denote by $P_x : W^{s, p}(M_0) \to
W^{s-m, p}(M_0)$ the induced family of pseudodifferential
operators, and we note that it will still be a smooth family that
is uniformly bounded in norm. Note that it is possible to extend
$P_x$ to an operator on $M_0$ because its distribution kernel has
compact support.

Then choose the sequence of points $\{x_j\} \subset M_0$ and a
partition of unity $\phi_j \in \CIc(M_0)$ as in Lemma
\ref{lemma.Shubin}. In particular, $\phi_j$ will have support in
$B(x_j, 2\epsilon)$. Also, let $\psi_j : B(x_j,4\epsilon) \to
B_{\RR^n}(0,4\epsilon)$ denote the normal system of coordinates
induced by the exponential maps $\exp_x : T_x M_0 \to M_0$ and
some fixed isometries $T_xM_0 \simeq \RR^n$. Then all derivatives
of $\psi_j \circ \psi_k^{-1}$ are bounded on their domain of
definition, with a bound that may depend on $\epsilon$ but does
not depend on $j$ and $k$ \cite{CGT, Shubin}.

Let
\begin{equation*}
    \nu_{s,p}(u)^p := \sum_j \|(\phi_j u) \circ
    \psi_j^{-1}\|_{W^{s,p}(\RR^n)}^p.
\end{equation*}
be one of the several equivalent norms defining the topology on
$W^{s, p}(M_0)$ (see Proposition \ref{prop.alt.desc1} and Equation
\eqref{sreal}. It is enough to prove that
\begin{multline}\label{mult.cont}
    \nu_{s,p}(Pu)^p := \sum_j \|(\phi_j Pu) \circ
    \psi_j^{-1}\|_{W^{s,p}(\RR^n)}^p \\ \le C \sum_j \|(\phi_j u)
    \circ \psi_j^{-1}\|_{W^{s,p}(\RR^n)}^p =: C\nu_{s,p}(u)^p,
\end{multline}
for some constant $C$ independent of $u$.

We now prove this statement. Indeed, for the reasons explained
below, we have the following inequalities.
\begin{multline*}
    \sum_j \|(\phi_j Pu) \circ \psi_j^{-1}\|_{W^{s,p}(\RR^n)}^p \le
    C \sum_{j, k} \|(\phi_j P \phi_k u ) \circ
    \psi_j^{-1}\|_{W^{s,p}(\RR^n)}^p \\ = C \sum_{j, k} \|(\phi_j
    P_{x_j} \phi_k u ) \circ \psi_j^{-1}\|_{W^{s,p}(\RR^n)}^p \le C
    \sum_{j, k} \|(\phi_j \phi_k u ) \circ
    \psi_j^{-1}\|_{W^{s,p}(\RR^n)}^p \\ \le C \sum_{j} \|(\phi_j u )
    \circ \psi_j^{-1}\|_{W^{s,p}(\RR^n)}^p = C \nu_{s,p}(u)^p.
\end{multline*}
Above, the first and last inequalities are due to the fact that
the family $\phi_j$ is uniformly locally finite, that is, there
exists a constant $\kappa$ such that at any given point $x$, at
most $\kappa$ of the functions $\phi_j(x)$ are different from
zero. The first equality is due to the support assumptions on
$\phi_j$, $\phi_k$, and $P_{x_j}$. Finally, the second inequality
is due to the fact that the operators $P_{x_j}$ are continuous,
with norms bounded by a constant independent of $j$, as explained
above. We have therefore proved that $P = q(a) \in \Psi_{1, 0,
\VV}^m(M_0)$ defines a bounded operator $W^{s, p}(M_0) \to W^{s-m,
p}(M_0)$, provided that the Schwartz kernel of $P$ has support in
a set of the $(M_0^2)_{\epsilon}$, for $\epsilon <
r_{\mathrm{inj}}(M_0)/9$.

Assume now that $P \in \Psi_{\VV}^{-\infty}(M_{0})$. We shall
check that $P$ is bounded as a map $W^{2k, p}(M_0) \to W^{-2k,
p}(M_0)$. For $k = 0$, this follows from the fact that the
Schwartz kernel of $P$ is given by a smooth function $k(x, y)$
such that $\int_{M_0}|k(x, y)|d\vol_g(x)$ and $\int_{M_0}|k(x,
y)|d\vol_g(y)$ are uniformly bounded in $x$ and $y$. For the other
values of $k$, it is enough to prove that the bilinear form
\begin{equation*}
    W^{2k, p}(M_0) \times W^{2k, p}(M_0) \ni (u, v) \mapsto \<Pu, v\> \in
    \CC
\end{equation*}
is continuous. Choose $Q$ a parametrix of $\Delta^k$ and let $R =
1 - Q\Delta^k$ be as above. Let $R' = 1 - \Delta^kQ \in
\Psi_{\VV}^{-\infty}(M_{0}).$ Then
\begin{equation*}
    \<Pu, v\> = \<(Q P Q)\Delta^{k} u, \Delta^{k}v\> + \<(Q P R) u,
    \Delta^{k}v\> + \<(R' P Q)\Delta^{k} u, v\> + \<(R' P R) u, v\>,
\end{equation*}
which is continuous since $QPQ, Q P R, R' P Q$, and $R' P R$ are
in $\Psi_{\VV}^{-\infty}(M_{0})$ and hence they are continuous
on $L^p(M_0)$ and because $\Delta^k : W^{2k, p}(M_0) \to L^p(M_0)$
is continuous.

Since any $P \in \Psi_{1, 0, \VV}^{m}(M_{0})$ can be written $P =
P_1 + P_2$ with $P_2 \in \Psi_{\VV}^{-\infty}(M_{0})$ and $P_1 =
q(a) \in \Psi_{1, 0, \VV}^{m}(M_{0})$ with support arbitrarily
close to the diagonal in $M_0$, the result follows.
\end{proof}

We obtain the following standard description of Sobolev spaces.

\begin{theorem}\label{theorem.three}\
Let $s \in \RR_+$ and $p \in (1, \infty)$. We have that $u \in
W^{s, p}(M_0)$ if, and only if, $u \in L^p(M_0)$ and $Pu \in
L^p(M_0)$ for any $P \in \Psi_{1, 0, \VV}^s(M_0)$. The norm $u \mapsto
\|u\|_{L^p(M_0)} + \|Pu\|_{L^p(M_0)}$ is equivalent to the
original norm on $W^{s,p}(M_0)$ for any elliptic $P \in \Psi_{1,
0, \VV}^s(M_0)$.

Similarly, the map $T: L^p(M_0) \oplus L^p(M_0) \ni (u, v) \mapsto u + Pv
\in W^{-s, p}(M_0)$ is surjective and identifies $W^{-s, p}(M_0)$ with
the quotient $(L^p(M_0) \oplus L^p(M_0))/\ker(T)$.
\end{theorem}

\begin{proof}\
Clearly, if $u \in W^{s, p}(M_0)$, then $Pu, u \in L^p(M_0)$. Let
us prove the converse. Assume $Pu, u \in L^p(M_0)$. Let $Q \in
\Psi_{1, 0, \VV}^{-s}(M_0)$ be a parametrix of $P$ and let $R, R'
\in \Psi_{\VV}^{-\infty}(M_0)$ be defined by $R := 1 - QP$ and $R'
= 1 - PQ$.  Then $u = QPu + Ru$. Since both $Q, R : L^p(M_0) \to
W^{s, p}(M_0)$ are defined and bounded, $u \in W^{s, p}(M_0)$ and
$\|u\|_{W^{s, p}(M_0)} \le C \big( \|u\|_{L^p(M_0)} +
\|Pu\|_{L^p(M_0)} \big )$. This proves the first part.

To prove the second part, we observe that the mapping
\begin{equation*}
    W^{s, q}(M_0) \ni u \mapsto (u, Pu) \in L^q(M_0)
    \oplus L^q(M_0), \qquad q^{-1} + p^{-1} = 1,
\end{equation*}
is an isomorphism onto its image. The result then follows by
duality using also the Hahn-Banach theorem.
\end{proof}

We conclude our paper with the sketch of a regularity results for
solutions of elliptic equations. Recall the Sobolev spaces with
weights $\rho^sW^{s, p}(\Omega_0)$ introduced in Equation
\eqref{eq.def.weights}.

\begin{theorem}\label{thm.reg2}\
Let $P \in \DiffV{m}(M)$ be an order $m$ elliptic operator on
$M_0$ generated by $\VV$. Let $u \in \rho^{s}W^{r, p}(M_0)$ be
such that $Pu \in \rho^{s}W^{t, p}(M_0)$, $s, r, t \in \RR$, $1 <
p < \infty$. Then $u \in \rho^{s}W^{t + m, p}(M_0)$.
\end{theorem}

\begin{proof}\
Let $Q \in \Psi_{\VV}^{-\infty}(M_0)$ be a parametrix of $P$. Then
$R = I - QP \in \Psi_{\VV}^{-\infty}(M_0)$. This gives $u = Q(Pu)
+ Ru$. But $Q(Pu) \in \rho^{s}W^{t + m, p}(M_0)$, by Theorem
\ref{thm.Wkp.cont}, because $Pu \in \rho^{s}W^{t, p}(M_0)$.
Similarly, $Ru \in \rho^{s}W^{t + m, p}(M_0)$. This completes the
proof.
\end{proof}

Note that the above theorem was already proved in the case $t \in
\ZZ$ and $m = 2$, using more elementary methods, as part of
Theorem \ref{thm.reg1}. The proof here is much shorter, however,
it attests to the power of pseudodifferential operator algebra
techniques.

\def\cprime{$'$} \def\cprime{$'$} \def\cprime{$'$} \def\cprime{$'$}
  \def\cprime{$'$}

\end{document}